\UseRawInputEncoding
\documentclass[10pt]{article}
\usepackage{bm,latexsym,amsmath,amsthm,amssymb}
\usepackage{graphicx,indentfirst}
\usepackage[colorlinks=true]{hyperref}
\newtheorem{theorem}{Theorem}[section]
\newtheorem{corollary}{Corollary}

\newtheorem{lemma}[theorem]{Lemma}

\theoremstyle{definition}

\newtheorem{remark}{Remark}

\begin{document}
\title{\bf A two-species competition model with mixed dispersal and free boundaries in time-periodic environment}
\author{ Qiaoling Chen\thanks{Corresponding author.\newline
\mbox{}\qquad E-mail: qiaolingf@126.com}
\\
\small School of Mathematics and Information Science, Shaanxi Normal University,
Xi'an 710062, PR China\\
\small School of Science, Xi'an Polytechnic University, Xi'an 710048, PR China\\
Fengquan Li\\
\small School of Mathematical Sciences, Dalian University of Technology,
Dalian 116024, PR China\\
Sanyi Tang\\
\small School of Mathematics and Information Science, Shaanxi Normal University,
Xi'an 710062, PR China\\
Feng Wang\\
\small School of Mathematics and Statistics, Xidian University
Xi'an 710071, PR China
}
\date{}
\maketitle \baselineskip 5pt
\begin{center}
\begin{minipage}{130mm}
{{\bf Abstract.} This paper is concerned with a Lotka-Volterra type competition model with free boundaries in time-periodic environment. One species is assumed to adopt nonlocal dispersal and the other
one adopts mixed dispersal, which is a combination of both random dispersal and nonlocal dispersal. We show that this free boundary problem with more general growth functions admits a unique solution defined for all time.  A spreading-vanishing dichotomy is obtained and criteria for spreading and vanishing are provided. Moreover, under the weak competition condition we provide the long-time asymptotic behavior of solution when spreading occurs.

}

\end{minipage}
\end{center}

\baselineskip=15pt

\section{Introduction}

In this paper, we study the dynamical behavior of the solution
$(u(t, x), v(t, x), g(t), h(t))$ to the following Lotka-Volterra type competition model with mixed dispersal and free boundaries in time-periodic environment
\begin{align*}\label{1.1}
\left\{\begin{array}{l}
\partial_{t}u=d_{1}\left(\int_{g(t)}^{h(t)}J(x-y)u(t,y)dy-u\right)+u(a(t)-u-b(t)v),\\[5pt]
\qquad\qquad\qquad\qquad\qquad\qquad\qquad\qquad\qquad\qquad\qquad\qquad
t>0,~g(t)<x<h(t),\\[5pt]
\partial_{t}v=d_{2}\left[\tau \partial_{x}^{2}v+(1-\tau)\left(\int_{g(t)}^{h(t)}J(x-y)v(t,y)dy-v\right)\right]
+v(c(t)-v-d(t)u),\\[5pt]
\qquad\qquad\qquad\qquad\qquad\qquad\qquad\qquad\qquad\qquad\qquad\qquad
t>0,~ g(t)<x<h(t),\\[5pt]
u(t, g(t))=u(t,h(t))=v(t, g(t))=v(t,h(t))=0, \quad t\geq0,\\[5pt]
h'(t)=-\mu v_{x}(t, h(t))
+\rho_{1}\int_{g(t)}^{h(t)}\int_{h(t)}^{+\infty}J(x-y)u(t,x)dydx\\[5pt]
\qquad\quad
+\rho_{2}\int_{g(t)}^{h(t)}\int_{h(t)}^{+\infty}J(x-y)v(t,x)dydx, \quad t\geq0,\\[5pt]
g'(t)=-\mu v_{x}(t, g(t))
-\rho_{1}\int_{g(t)}^{h(t)}\int_{-\infty}^{g(t)}J(x-y)u(t,x)dydx\\[5pt]
\qquad\quad
-\rho_{2}\int_{g(t)}^{h(t)}\int_{-\infty}^{g(t)}J(x-y)v(t,x)dydx, \quad t\geq0,\\[5pt]
u(0,x)=u_0(x),~v(0, x)=v_0(x), \quad |x|\leq h_{0},\\[5pt]
h(0)=-g(0)=h_{0}.
\end{array}\right.
\tag{1.1}
\end{align*}
Here $u(t,x)$ and $v(t,x)$ represent the population densities of two competing species; the positive constants $d_1, d_2$ are dispersal rates of $u, v$ and the constant $0<\tau\leq 1$ measures the fraction of individuals adopting random dispersal; $h_{0}$, $\mu$ and $\rho_{1}$ are positive constants, $\rho_{2}$ is a nonnegative constant, $\rho_{2}>0$ when $\tau<1$ and $\rho_{2}=0$ when $\tau=1$;
the kernel function $J: \mathbb{R}\rightarrow \mathbb{R}$ satisfies that
\begin{align*}\label{J}
\tag{J}
\begin{array}{c}
J~\mbox{is~Lipschitz~continuous},~
J(x)\geq 0,~J(0)>0,~\int_{\mathbb{R}}J(x)dx=1, \\[3pt]
J~\mbox{is~symmetric~and}~\sup_{\mathbb{R}}J<\infty;
\end{array}
\end{align*}
$a(t),c(t)$ represent the intrinsic growth rates of species, $b(t),d(t)$ represent competition between species and they satisfy that
\begin{align*}
a(t),b(t),c(t),d(t)~\mbox{are~positive}~T\mbox{-periodic~functions}\\
\mbox{and}~a,b\in C([0,T]),~c,d\in C^{\frac{\alpha}{2}}([0,T])~\mbox{for}~0<\alpha<1;
\end{align*}
the initial functions $u_0$ and $v_0$ satisfy
\begin{align*}\label{1.2}
\left\{\begin{array}{l}
u_0\in C^{1-}([-h_{0},h_{0}]),\quad u_0(\pm h_0)=0,\quad u_0>0 \quad \mbox{in}~ (-h_{0}, h_0),\\[3pt]
v_0\in C^{2}([-h_{0},h_{0}]), \quad v_0(\pm h_0)=0, \quad v_0>0 \quad \mbox{in}~ (-h_{0}, h_0),
\end{array}\right.
 \tag{1.2}
\end{align*}
where $C^{1-}([-h_{0},h_{0}])$ is defined as the Lipschitz continuous function space.

Ecologically, problem $(\ref{1.1})$ describes the dynamical process of two competing species which spread and invade to new environment with daily or seasonal changes via the same free boundaries. All the individuals in the population $u$ adopt nonlocal dispersal, while in the population $v$ a fraction of individuals adopt nonlocal dispersal and the remaining fraction assumes random dispersal. The latter strategy is called mixed dispersal, which was first proposed by Kao et al. \cite{kls12}.
We assume that the spreading fronts expand at a speed that is proportional to the outward flux of the population of the two species at the front, which give rise to the free boundary conditions in $(\ref{1.1})$.
Problem $(\ref{1.1})$ is a variation of the following two species competition system studied in \cite{kls12}:
\begin{align*}
\left\{\begin{array}{l}
\partial_{t}u=d_{1}\left(\int_{\mathbb{R}^{N}}J(x-y)u(t,y)dy-u\right)
+u(a(x)-u-v), \\[5pt]
\partial_{t}v=d_{2}\left[\tau \partial_{x}^{2}v+(1-\tau)\left(\int_{\mathbb{R}^{N}}J(x-y)v(t,y)dy-v\right)\right]
+v(a(x)-u-v).
\end{array}\right.
\end{align*}
They investigated how the mixed dispersal affects the invasion of a single species and how the mixed dispersal strategies will evolve in spatially periodic but temporally constant environment. A complete classification of the global dynamics of two-species competition mode with mixed dispersals was studied in
\cite{bl18}. If both $u,v$ adopt random dispersal, the existence and stability of time periodic traveling waves to $(\ref{1.1})$ with $x\in \mathbb{R}$ have been studied in \cite{bw13,zr11}.

If $\tau=0$ and $a(t),b(t),c(t),d(t)$ are constants, $(\ref{1.1})$ reduces to a two species nonlocal diffusion system with free boundaries studied by Du et al. \cite{dwz19}. They proved the model has a unique global solution, established a spreading-vanishing dichotomy and obtained criteria for spreading and vanishing. Moreover, for the weak competition case they determined the long-time asymptotic limit of the solution when spreading happens. If $\tau=1$ and $a(t),b(t),c(t),d(t)$ are constants, $(\ref{1.1})$ becomes a free boundary problem of ecological model with nonlocal and local diffusions considered in \cite{waw181,waw182}. They also obtained well-posedness of solutions and spreading-vanishing results. Moreover, Cao et al. \cite{clwan19} considered a nonlocal diffusion Lotka-Volterra type competition model with free boundaries in the homogeneous environment, which consists of a native species distributing in the whole space $\mathbb{R}$ and an invasive species. Some free boundary problems for epidemic models with nonlocal dispersals have been recently studied in \cite{dn20,zzld20}.

In the absence of the species $v$ (i.e., $v\equiv 0$) and $a(t)$ is a constant, $(\ref{1.1})$ reduces to the following nonlocal dispersal model with free boundaries
\begin{align*}\label{1.3}
\left\{\begin{array}{l}
\partial_{t}u=d_{1}\left(\int_{g(t)}^{h(t)}J(x-y)u(t,y)dy-u\right)+u(a-u),\quad t>0,~g(t)<x<h(t),\\[3pt]
u(t, g(t))=u(t,h(t))=0, \quad t\geq0,\\[5pt]
h'(t)=
\rho_{1}\int_{g(t)}^{h(t)}\int_{h(t)}^{\infty}J(x-y)u(t,x)dydx, \quad t\geq0,\\[5pt]
g'(t)=
-\rho_{1}\int_{g(t)}^{h(t)}\int_{-\infty}^{g(t)}J(x-y)u(t,x)dydx, \quad t\geq0,\\[5pt]
u(0,x)=u_0(x),\quad |x|\leq h_{0},\\[5pt]
h(0)=-g(0)=h_{0}.
\end{array}\right.
 \tag{1.3}
\end{align*}
which has been studied in \cite{cdll19}. Problem $(\ref{1.3})$ is a nature extension of the random dispersal model with free boundary in \cite{dl10}, and similar results including the existence and uniqueness of global solutions for more general growth function $f(t,x,u)$ and the spreading-vanishing results in the homogeneous environment were obtained in \cite{cdll19}, from which one can see that the nonlocal dispersal brings many essential difficulties in analysis. The spreading speed of $(\ref{1.3})$ was determined in \cite{dlz19} when spreading happens. After this paper is completed, we learned of the paper \cite{zlz20}, where $(\ref{1.3})$ with the assumptions $a(t,x)=\alpha(t)+\beta(x)$ ($\alpha(t)$ is $T$-periodic) and $\text{supp}~ J\subset[-r_{0},r_{0}]$ was studied.

Since the work of Du and Lin \cite{dl10}, the random dispersal models with free boundary(ies) have been studied extensively. For example, the model in \cite{dl10} has been extended to free boundary problems for single species models in heterogeneous environment and time-periodic environment, or with general nonlinear term, advection term and time delay, we refer the readers to \cite{bdk12,dpw17,dgp13,dbl15,gll15,k14,llz14,llo15,rz19,sf19,w143,chw18,zx14}
and references therein. Moreover, two-species Lotka-Volterra type competition problems and predator-prey problems with free boundary(ies) have been considered in the homogeneous environment or heterogeneous time-periodic environment, e.g., \cite{clw17,dl13,dwz17,gw12,gw15,tir18,w14,waz16,wz18,wz15,zzl17}. Free boundary problems for epidemic models with random dispersal \cite{clwy17,gklz15,lzh17} and time delay \cite{cltw21} have also been considered recently.

The aim of this paper is to study the well-posedness and long-time behavior of solution to problem $(\ref{1.1})$. We first investigate the existence and uniqueness of solutions to $(\ref{1.1})$ with more general growth functions. To achieve it, we shall establish the
maximum principle for linear parabolic equations with mixed dispersal, and prove that the nonlinear parabolic equations with mixed dispersal (see $(\ref{2.5})$) admit a unique positive solution under the assumption that $g^{\prime}(t),h^{\prime}(t)$ and $u(t,x)$ are only continuous functions by approximation method, which plays an important role in the process of using the fixed point theorem (see Lemma $\ref{l2.5}$).
Then we establish a spreading-vanishing dichotomy and criteria for spreading and vanishing. To discuss the spreading and vanishing, we need to consider the existence and properties of principle eigenvalue of
time-periodic parabolic-type eigenvalue problems with random/mixed dispersal. Since the intrinsic growth rates $a(t)$ and $c(t)$ are independent of spatial variable, we can transform the parabolic-type eigenvalue problems into elliptic-type eigenvalue problems. This transformation is also used in discussing the asymptotic behavior of solution (see Theorem $\ref{t4.4}$). Moreover, by the comparison principle established in Lemma $\ref{l3-3}$, we discuss the asymptotic stability and uniqueness of $T$-periodic solutions to the nonlocal and mixed dispersal equations in $\mathbb{R}$ (Lemma $\ref{l4.6}$), which are used to determine the long-time asymptotic behavior of solution when spreading occurs under the weak competition condition (Theorem $\ref{t4.7}$).

The rest of the paper is organized as follows. In Section 2, we establish the global existence and uniqueness of solutions to problem $(\ref{1.1})$ with more general growth functions. The comparison principle in the moving domain and the discussions on eigenvalue problems are given in Section 3. In Section 4, we investigate spreading and vanishing of species.

\section{Well-posedness}

In this section, we give the global well-posedness of solutions to problem $(\ref{1.1})$ with more general growth functions. More precisely, we assume that $f_{i}(t,x,u,v)$ $(i=1,2)$ satisfy the following assumptions:

$(\textbf{f1})$ $f_{1}(t,x,0,v),f_{2}(t,x,u,0)\equiv0$, and there exists a constant $K>0$ such that
$f_{1}(t,x,u,v)<0$ for all $u>K$, $v\geq 0$ and $(t,x)\in \mathbb{R}^{+}\times \mathbb{R}$, and
$f_{2}(t,x,u,v)<0$ for all $u\geq 0$, $v>K$ and $(t,x)\in \mathbb{R}^{+}\times \mathbb{R}$;

$(\textbf{f2})$ For any given $T, l, K_{1}, K_{2}>0$, there exists a constant $L=L(T,l,K_{1},K_{2})$ such that
$$
\|f_{2}(\cdot,x,u,v)\|_{C^{\frac{\alpha}{2}}([0,T])}\leq L
$$
for all $x\in [-l,l]$, $u\in [0,K_{1}]$ and $v\in [0,K_{2}]$;

$(\textbf{f3})$ For any $K_{1},K_{2}>0$, there exists a constant $L^{*}=L^{*}(K_{1},K_{2})>0$ such that
$$
|f_{i}(t,x,u,v)-f_{i}(t,y,u,v)|\leq L^{*}|x-y|
$$
for all $u\in [0,K_{1}]$, $v\in [0,K_{2}]$ and all $(t,x,y)\in \mathbb{R}^{+}\times \mathbb{R}\times \mathbb{R}$;

$(\textbf{f4})$ $f_{i}(t,x,u,v)$ is locally Lipschitz in $u,v\in \mathbb{R}^{+}$
uniformly for $(t,x)\in \mathbb{R}^{+}\times \mathbb{R}$, i.e.,
for any $K_{1},K_{2}>0$, there exists a constant $\hat{L}=\hat{L}(K_{1},K_{2})>0$ such that
$$
|f_{i}(t,x,u_{1},v_{1})-f_{i}(t,x,u_{2},v_{2})|\leq \hat{L}(|u_{1}-u_{2}|+|v_{1}-v_{2}|)
$$
for all $u_{1},u_{2}\in [0,K_{1}]$, $v_{1},v_{2}\in [0,K_{2}]$ and all $(t,x)\in \mathbb{R}^{+}\times \mathbb{R}$.

It is easy to check that the growth functions in $(\ref{1.1})$ satisfy the conditions $(\textbf{f1})-(\textbf{f4})$. We consider the following free boundary problem
\begin{align*}\label{2.1}
\left\{\begin{array}{l}
\partial_{t}u=d_{1}\left(\int_{g(t)}^{h(t)}J(x-y)u(t,y)dy-u\right)+f_{1}(t,x,u,v), \quad t>0,~g(t)<x<h(t),\\[5pt]
\partial_{t}v=d_{2}\left[\tau \partial_{x}^{2}v+(1-\tau)\left(\int_{g(t)}^{h(t)}J(x-y)v(t,y)dy-v\right)\right]\\[5pt]
\qquad\quad+f_{2}(t,x,u,v),
\quad t>0,~ g(t)<x<h(t),\\[5pt]
u(t, g(t))=u(t,h(t))=v(t, g(t))=v(t,h(t))=0, \quad t\geq0,\\[5pt]
h'(t)=-\mu v_{x}(t, h(t))
+\rho_{1}\int_{g(t)}^{h(t)}\int_{h(t)}^{\infty}J(x-y)u(t,x)dydx\\[5pt]
\qquad\quad+\rho_{2}\int_{g(t)}^{h(t)}\int_{h(t)}^{\infty}J(x-y)v(t,x)dydx, \quad t\geq0,\\[5pt]
g'(t)=-\mu v_{x}(t, g(t))
-\rho_{1}\int_{g(t)}^{h(t)}\int_{-\infty}^{g(t)}J(x-y)u(t,x)dydx\\[5pt]
\qquad\quad-\rho_{2}\int_{g(t)}^{h(t)}\int_{-\infty}^{g(t)}J(x-y)v(t,x)dydx, \quad t\geq0,\\[5pt]
u(0,x)=u_0(x), v(0, x)=v_0(x), \quad |x|\leq h_{0},\\[5pt]
h(0)=-g(0)=h_{0},
\end{array}\right.
\tag{2.1}
\end{align*}
Throughout the paper, we denote
$\Omega_{T_{0}}^{g,h}=(0,T_{0}]\times (g(t), h(t))$,
$D_{T_{0}}=(0,T_{0}]\times (-1, 1)$ and $a_{T}=\frac{1}{T}\int_{0}^{T}a(t)dt$. Under the transform
$x(t,z)=\frac{(h(t)-g(t))z+h(t)+g(t)}{2}$, we always denote
$\tilde{f}(t,z)=f(t,x(t,z))=f(t,\frac{(h(t)-g(t))z+h(t)+g(t)}{2})$. $C^{1,1-}(\overline{\Omega}_{T_{0}}^{g,h})$ denotes the class of functions that
are $C^{1}$ in $t$ and Lipschitz continuous in $x$. The main result of this section is stated in the following theorem.

\begin{theorem}\label{t2.1} Assume that $(\ref{J})$ and (\textbf{f1})-(\textbf{f4}) hold.
For any given $(u_{0},v_{0})$ satisfying $(\ref{1.2})$, the problem $(\ref{2.1})$ admits a unique global solution
$(u, v, g, h)$ defined on $[0,T_{0}]$ for any $0<T_{0}<\infty$ and
\begin{align*}\label{2.2}
\begin{array}{rl}
&(u, v, g, h)\in C^{1,1-}(\overline{\Omega}_{T_{0}}^{g,h})\times C^{1+\frac{\alpha}{2},2+\alpha}(\Omega_{T_{0}}^{g,h})\times [C^{1+\frac{\alpha}{2}}([0, T_{0}])]^{2},\\[3pt]
&0<u\leq K_{1},~
0<v\leq K_{2},
\quad\forall~
(t,x)\in\Omega_{T_{0}}^{g,h}, \\[3pt]
&0<-v_{x}(t,h(t)),~v_{x}(t,g(t))\leq K_{3},\quad 0<t\leq T_{0},
\end{array}
\tag{2.2}
\end{align*}
where
\begin{align*}
\begin{array}{c}
K_{1}:=\max\{\|u_{0}\|_{L^{\infty}}, K\},\quad K_{2}:=\max\{\|v_{0}\|_{L^{\infty}}, K\}, \\[3pt] K_{3}:=2K_{2}\max\left\{\sqrt{\frac{\hat{L}+d_{2}(1-\tau)}{2d_{2}\tau}},
\frac{4\|v_{0}\|_{C^{1}([-h_{0},h_{0}])}}{3K_{2}}\right\}
\end{array}
\end{align*}
and $\hat{L}=\hat{L}(K_{1},K_{2})$ is the Lipschitz constant defined in $(\textbf{f4})$.
\end{theorem}

To prove Theorem $\ref{t2.1}$, we first establish the maximum principle for linear parabolic equations with mixed dispersal. For some $h_{0}, T_{0}$, we define
\begin{align*}
&\mathbb{H}_{T_{0}}^{h_{0}}
:=\{h\in C^{1}([0,T_{0}]): h(0)=h_{0},
~0<h^{\prime}(t)\leq R(t)\}, \\[3pt]
&\mathbb{G}_{T_{0}}^{h_{0}}
:=\{g\in C^{1}([0,T_{0}]):~-g\in \mathbb{H}_{T_{0}}^{h_{0}}\}
\end{align*}
with
\begin{align*}
R(t):=\mu K_{3}+2(h_{0}\rho_{1}K_{1}+h_{0}\rho_{2}K_{2}+\mu K_{3})e^{(\rho_{1}K_{1}+\rho_{2}K_{2})t}.
\end{align*}

\begin{lemma}\label{l2.2} (Maximum Principle) Assume that $(\ref{J})$ holds and
$(g,h)\in \mathbb{G}_{T_{0}}^{h_{0}}\times\mathbb{H}_{T_{0}}^{h_{0}}$. If $v(t,x)\in C^{1,2}(\Omega_{T_{0}}^{g,h})\cap C(\overline{\Omega}_{T_{0}}^{g,h})$ satisfies, for some $c\in L^{\infty}(\Omega_{T_{0}}^{g,h})$,
\begin{align*}\label{2.3}
\left\{\begin{array}{l}
\partial_{t}v\geq
d_{2}\left[\tau \partial_{x}^{2}v
+(1-\tau)\left(\int_{g(t)}^{h(t)}J(x-y)v(t,y)dy-v\right)\right]+c(t,x)v, \quad (t,x)\in\Omega_{T_{0}}^{g,h},\\[5pt]
v(t,g(t))\geq0,~v(t,h(t))\geq0,\quad t\in(0, T_{0}],\\[5pt]
v(0, x)\geq 0, \quad x\in[-h_{0}, h_{0}],
\end{array}\right.
 \tag{2.3}
\end{align*}
then $v(t,x)\geq 0$ for all $(t,x)\in\overline{\Omega}_{T_{0}}^{g,h}$.
Moreover, if $v(0,x)\not\equiv 0$ in $[-h_{0}, h_{0}]$, then
$v(t,x)>0$ in $\Omega_{T_{0}}^{g,h}$.
\end{lemma}

\begin{proof} $(i)$ Let $\omega(t,x)=e^{-kt}v(t,x)$, where
$k>0$ is a constant chosen large enough such that $-k+c(t,x)<0$ for all
$(t,x)\in \Omega_{T_{0}}^{g,h}$.
Then
\begin{align*}
\begin{array}{l}
\partial_{t}\omega\geq
d_{2}\left[\tau \partial_{x}^{2}\omega
+(1-\tau)\int_{g(t)}^{h(t)}J(x-y)\omega(t,y)dy\right]
+[-k-d_{2}(1-\tau)+c(t,x)]\omega.
\end{array}
\end{align*}
We are now in a position to prove that $\omega\geq 0$
in $\overline{\Omega}_{T_{0}}^{g,h}$.

Suppose that $\omega_{\inf}=\inf_{(t,x)\in \overline{\Omega}_{T_{0}}^{g,h}}\omega(t,x)<0$.
By $(\ref{2.3})$, $\omega\geq 0$ on the parabolic boundary of $\overline{\Omega}_{T_{0}}^{g,h}$,
and hence there exists $(t_{*}, x_{*})\in \Omega_{T_{0}}^{g,h}$ such that
$\omega_{\inf}=\omega(t_{*}, x_{*})<0$.
Since $\partial_{t}\omega(t_{*}, x_{*})\leq 0$,
$\partial_{x}^{2}\omega(t_{*}, x_{*})\geq 0$,
then
\begin{align*}
\begin{array}{rl}
\partial_{t}\omega(t_{*}, x_{*})
&\geq d_{2}\left[\tau \partial_{x}^{2}\omega(t_{*}, x_{*})
+(1-\tau)\int_{g(t_{*})}^{h(t_{*})}J(x_{*}-y)\omega(t_{*},y)dy\right]\\[5pt]
&\quad+[-k-d_{2}(1-\tau)+c(t_{*},x_{*})]\omega(t_{*}, x_{*})\\[5pt]
&\geq d_{2}\tau \partial_{x}^{2}\omega(t_{*}, x_{*})
+d_{2}(1-\tau)w_{\inf}\int_{\mathbb{R}}J(x_{*}-y)dy\\[5pt]
&\quad+[-k-d_{2}(1-\tau)+c(t_{*},x_{*})]w_{\inf}\\[5pt]
&=d_{2}\tau \partial_{x}^{2}\omega(t_{*}, x_{*})
+[-k+c(t_{*},x_{*})]w_{\inf}.
\end{array}
\end{align*}
Since $[-k+c(t_{*},x_{*})]\omega_{\inf}>0$,
we can get a contradiction. Thus,
$\omega(t,x)\geq 0$ in $\overline{\Omega}_{T_{0}}^{g,h}$, which implies that
\begin{align*}\label{2.4}
v(t,x)\geq 0\quad \mbox{for~all}~(t,x)\in\overline{\Omega}_{T_{0}}^{g,h}.
 \tag{2.4}
\end{align*}

$(ii)$ Now assume that $v(0,x)\not\equiv 0$ in $[-h_{0},h_{0}]$.
By $(\ref{2.4})$ and the fact $J(x)\geq 0$, we have
\begin{align*}
\begin{array}{rl}
\partial_{t}v
&\geq d_{2}\left[\tau \partial_{x}^{2}v
+(1-\tau)\left(\int_{g(t)}^{h(t)}J(x-y)v(t,y)dy-v\right)\right]+c(t,x)v\\[5pt]
&\geq d_{2}[\tau\partial_{x}^{2}v-(1-\tau)v]
+c(t,x)v
=d_{2}\tau\partial_{x}^{2}v
+[c(t,x)-d_{2}(1-\tau)]v.
\end{array}
\end{align*}
Define the transform
\begin{align*}
\begin{array}{rl}
x(t,z)=\frac{(h(t)-g(t))z+h(t)+g(t)}{2}, \quad \mbox{that~is,}\quad
z(t,x)=\frac{2x-g(t)-h(t)}{h(t)-g(t)},
\end{array}
\end{align*}
and let $\tilde{v}(t,z)=v(t,x(t,z))$ and $\tilde{c}(t,z)=c(t,x(t,z))$, then $\tilde{v}(t,z)$ satisfies
\begin{align*}
\left\{\begin{array}{l}
\partial_{t}\tilde{v}\geq
d_{2}\tau\xi(t)\partial_{z}^{2}\tilde{v}+\eta(t,z)\partial_{z}\tilde{v}
+[\tilde{c}(t,z)-d_{2}(1-\tau)]\tilde{v}, \quad (t,z)\in D_{T_{0}},\\[5pt]
\tilde{v}(t,-1)\geq0,~\tilde{v}(t,1)\geq0,\quad t\in(0, T_{0}],\\[5pt]
\tilde{v}(0, z)=v(0, h_{0}z)\geq 0, \quad z\in[-1, 1],
\end{array}\right.
\end{align*}
where
\begin{align*}
\begin{array}{rl}
\xi(t)=\frac{4}{(h(t)-g(t))^{2}},\quad
\eta(t,z)=\frac{h^{\prime}(t)+g^{\prime}(t)}{h(t)-g(t)}
+\frac{(h^{\prime}(t)-g^{\prime}(t))z}{h(t)-g(t)}.
\end{array}
\end{align*}
By the classical maximum principle for parabolic equation, we know
$\tilde{v}(t,z)> 0$, $\forall~(t,z)\in D_{T_{0}}$.
Thus, $v(t,x)>0$ in $\Omega_{T_{0}}^{g,h}$. This completes the proof.
\end{proof}

Next, we shall prove that nonlinear parabolic equations with mixed dispersal (see $(\ref{2.5})$) admit a unique positive strong solution for given continuous function $u(t,x)$ and $C^{1}$-functions $(g(t),h(t))\in \mathbb{G}_{T_{0}}^{h_{0}}\times\mathbb{H}_{T_{0}}^{h_{0}}$. The proof is divided into two lemmas. First, in Lemma $\ref{l2.3}$ we establish the existence and uniqueness of positive classical solutions by applying the upper-lower solutions method, under the assumption that $u(t,x)$ are H\"{o}lder continuous and
$(g(t),h(t))\in \widehat{\mathbb{G}}_{T_{0}}^{h_{0}}\times\widehat{\mathbb{H}}_{T_{0}}^{h_{0}}$ with \begin{align*}
&\widehat{\mathbb{H}}_{T_{0}}^{h_{0}}
:=\{h\in C^{1+\frac{\alpha}{2}}([0,T_{0}]): h(0)=h_{0},
~0<h^{\prime}(t)\leq R(t)\}, \\[3pt]
&\widehat{\mathbb{G}}_{T_{0}}^{h_{0}}
:=\{g\in C^{1+\frac{\alpha}{2}}([0,T_{0}]):~-g\in \widehat{\mathbb{H}}_{T_{0}}^{h_{0}}\}.
\end{align*}
Then, in Lemma $\ref{l2.4}$ we get the desired result by the approximation method.

\begin{lemma}\label{l2.3} Suppose that $(\ref{J})$ holds, $(g,h)\in \widehat{\mathbb{G}}_{T_{0}}^{h_{0}}\times\widehat{\mathbb{H}}_{T_{0}}^{h_{0}}$,
$u\in C^{\frac{\alpha}{2},\alpha}(\overline{\Omega}_{T_{0}}^{g,h})$,
$f_{2}$ satisfies (\textbf{f1})-(\textbf{f4})
and $v_{0}$ satisfies $(\ref{1.2})$. Then for any $T_{0}>0$, the following problem
\begin{align*}\label{2.5}
\left\{\begin{array}{l}
\partial_{t}v=
d_{2}\left[\tau \partial_{x}^{2}v
+(1-\tau)\left(\int_{g(t)}^{h(t)}J(x-y)v(t,y)dy-v\right)\right]\\[5pt]
\qquad\quad+f_{2}(t,x,u,v), \quad
(t,x)\in\Omega_{T_{0}}^{g,h},\\[5pt]
v(t,g(t))=v(t,h(t))=0,\quad t\in(0, T_{0}],\\[5pt]
v(0, x)=v_{0}(x), \quad x\in[-h_{0}, h_{0}]
\end{array}\right.
 \tag{2.5}
\end{align*}
admits a unique solution $v(t,x)\in C^{1+\frac{\alpha}{2},2+\alpha}(\Omega_{T_{0}}^{g,h})$.
Moreover, $v(t,x)$ satisfies
\begin{align*}\label{2.6}
\begin{array}{rl}
&0<v(t,x)\leq K_{2} \quad\mbox{for}~(t,x)\in\Omega_{T_{0}}^{g,h},\\[5pt]
&0<-v_{x}(t,h(t)), v_{x}(t,g(t))\leq K_{3}
\quad\mbox{for}~t\in(0,T_{0}].
\end{array}
 \tag{2.6}
\end{align*}
\end{lemma}

\begin{proof}
For the existence and uniqueness, we mainly adopt the classical upper-lower solutions method. Since the mixed dispersal is considered, we give some details of the proof. A function $\bar{v}$ is called an upper solution of $(\ref{2.5})$ if
$\bar{v}\in C^{1,2}(\Omega_{T_{0}}^{g,h})\cap C(\overline{\Omega}_{T_{0}}^{g,h})$ satisfies
\begin{align*}
\left\{\begin{array}{l}
\partial_{t}\bar{v}\geq
d_{2}\left[\tau \partial_{x}^{2}\bar{v}
+(1-\tau)\left(\int_{g(t)}^{h(t)}J(x-y)\bar{v}(t,y)dy-\bar{v}\right)\right]\\[5pt]
\qquad\quad+f_{2}(t,x,u,\bar{v}),
\quad (t,x)\in\Omega_{T_{0}}^{g,h},\\[5pt]
\bar{v}(t,g(t))\geq 0,~\bar{v}(t,h(t))\geq0,\quad t\in(0, T_{0}],\\[5pt]
\bar{v}(0, x)\geq v_{0}(x), \quad x\in[-h_{0}, h_{0}],
\end{array}\right.
\end{align*}
and a function $\underline{v}$ is called a lower solution of $(\ref{2.5})$ if reversing all the above inequalities.

\emph{Step 1.}
We claim that, if $\bar{v}, \underline{v}$
are respectively nonnegative upper and lower solutions of $(\ref{2.5})$, then $(\ref{2.5})$ has a unique solution $v(t,x)$ satisfying
$\underline{v}(t,x)\leq v(t,x)\leq \bar{v}(t,x)$, $\forall(t,x)\in \overline{\Omega}_{T_{0}}^{g,h}$.

Indeed, since $u\in C^{\frac{\alpha}{2},\alpha}(\overline{\Omega}_{T_{0}}^{g,h})$ and $\bar{v}, \underline{v}\in C(\overline{\Omega}_{T_{0}}^{g,h})$, there exists a constant $M>0$ such that
$0\leq u,\bar{v}, \underline{v}\leq M$ for $(t,x)\in \overline{\Omega}_{T_{0}}^{g,h}$.
By (\textbf{f4}),
we have, for some constant $k>d_{2}(1-\tau)$,
\begin{align*}
|f_{2}(t,x,u,v_{1})-f_{2}(t,x,u,v_{2})|\leq [k-d_{2}(1-\tau)]|v_{1}-v_{2}|
\end{align*}
for any $(t,x)\in \overline{\Omega}_{T_{0}}^{g,h}$ and
$u, v_{1}, v_{2}\in [0,M]$.

For any $\vartheta\in C(\overline{\Omega}_{T_{0}}^{g,h})$ satisfying $\vartheta\in[0,M]$, we define a mapping $\Phi$ by $v=\Phi \vartheta$, where $v\in C^{\frac{1+\alpha}{2},1+\alpha}(\overline{\Omega}_{T_{0}}^{g,h})$ is the unique solution of
\begin{align*}\label{2.7}
\left\{\begin{array}{l}
\partial_{t}v-
d_{2}\tau \partial_{x}^{2}v+kv
=d_{2}(1-\tau)\left(\int_{g(t)}^{h(t)}J(x-y)\vartheta(t,y)dy-\vartheta\right)\\[5pt]
\qquad\qquad\qquad\qquad\quad
+f_{2}(t,x,u,\vartheta)+k\vartheta,
\qquad
(t,x)\in\Omega_{T_{0}}^{g,h},\\[5pt]
v(t,g(t))=v(t,h(t))=0,\quad t\in(0, T_{0}],\\[5pt]
v(0, x)=v_{0}(x), \quad x\in[-h_{0}, h_{0}].
\end{array}\right.
 \tag{2.7}
\end{align*}
The existence and uniqueness of $v\in C^{\frac{1+\alpha}{2},1+\alpha}(\overline{\Omega}_{T_{0}}^{g,h})$ is guaranteed by the $L^{p}$ theory for linear parabolic equation and the Sobolev imbedding theorem.
More precisely, let $\tilde{v}(t,z)=v(t,x(t,z))$, $\tilde{u}(t,z)=u(t,x(t,z))$, $\tilde{\vartheta}(t,z)=\vartheta(t,x(t,z))$ and $\tilde{f}_{2}(t,z,\tilde{u},\tilde{\vartheta})=f_{2}(t,x(t,z),\tilde{u},\tilde{\vartheta})$, then $(\ref{2.7})$ becomes
\begin{align*}\label{2.8}
\left\{\begin{array}{l}
\partial_{t}\tilde{v}-
d_{2}\tau\xi(t)\partial_{z}^{2}\tilde{v}-\eta(t,z)\partial_{z}\tilde{v}+k\tilde{v}\\[5pt]
=d_{2}(1-\tau)\left(\frac{h(t)-g(t)}{2}\int_{-1}^{1}J(\frac{h(t)-g(t)}{2}(z-s))
\tilde{\vartheta}ds-\tilde{\vartheta}\right)\\[5pt]
\quad+\tilde{f}_{2}(t,z,\tilde{u},\tilde{\vartheta})+k\tilde{\vartheta},
\qquad
(t,z)\in D_{T_{0}},\\[5pt]
\tilde{v}(t,-1)=\tilde{v}(t,1)=0,\quad t\in(0, T_{0}],\\[5pt]
\tilde{v}(0, z)=v_{0}(h_{0}z), \quad z\in[-1, 1].
\end{array}\right.
 \tag{2.8}
\end{align*}
Note that the right hand of the equation in $(\ref{2.8})$ is continuous in $\overline{D}_{T_{0}}$ and then belongs to $L^{p}(D_{T_{0}})$ with any $p>3$, $\xi(t)\in C([0,T_{0}])$ with $\|\xi\|_{L^{\infty}((0,T_{0}))}\leq \frac{1}{h_{0}^{2}}$
and $\|\eta\|_{L^{\infty}((0,T_{0}))}\leq \frac{2R(T_{0})}{h_{0}}$. Applying the $L^{p}$ theory to $(\ref{2.8})$ and the Sobolev imbedding theorem, we can obtain a unique solution $\tilde{v}\in W^{1,2}_{p}(D_{T_{0}})\hookrightarrow C^{\frac{1+\alpha}{2},1+\alpha}(\overline{D}_{T_{0}})$, and then get a unique solution $v\in C^{\frac{1+\alpha}{2},1+\alpha}(\overline{\Omega}_{T_{0}}^{g,h})$ to $(\ref{2.7})$.

We shall show that $\Phi$ is monotone in the sense that if any $\vartheta_{1}, \vartheta_{2}\in C(\overline{\Omega}_{T_{0}}^{g,h})$ satisfy
$0\leq \vartheta_{1},\vartheta_{2}\leq M$ and $\vartheta_{2}\geq \vartheta_{1}$, then
$\Phi \vartheta_{2}\geq \Phi \vartheta_{1}$.
To see that, let $w=\Phi \vartheta_{2}-\Phi \vartheta_{1}$, then $w$ satisfies
\begin{align*}\label{2.9}
\left\{\begin{array}{l}
\partial_{t}w-d_{2}\tau \partial_{x}^{2}w+kw=d_{2}(1-\tau)\left(\int_{g(t)}^{h(t)}J(x-y)(\vartheta_{2}(t,y)-\vartheta_{1}(t,y))dy
-(\vartheta_{2}-\vartheta_{1})\right)\\[5pt]
\qquad\qquad\qquad\qquad\quad
+f_{2}(t,x,u,\vartheta_{2})-f_{2}(t,x,u,\vartheta_{1})
+k(\vartheta_{2}-\vartheta_{1}), \quad (t,x)\in\Omega_{T_{0}}^{g,h},\\[5pt]
w(t,g(t))=w(t,h(t))=0,\quad t\in(0, T_{0}],\\[5pt]
w(0, x)=0, \quad x\in[-h_{0}, h_{0}].
\end{array}\right.
 \tag{2.9}
\end{align*}
Since the equation in $(\ref{2.9})$ satisfies
\begin{align*}
\begin{array}{rl}
&\partial_{t}w-
d_{2}\tau \partial_{x}^{2}w+kw\\[5pt]
&=d_{2}(1-\tau)\left(\int_{g(t)}^{h(t)}J(x-y)(\vartheta_{2}(t,y)-\vartheta_{1}(t,y))dy
-(\vartheta_{2}-\vartheta_{1})\right)\\[5pt]
&\quad+f_{2}(t,x,u,\vartheta_{2})-f_{2}(t,x,u,\vartheta_{1})
+k(\vartheta_{2}-\vartheta_{1})\\[5pt]
&\geq -d_{2}(1-\tau)(\vartheta_{2}-\vartheta_{1})
+f_{2}(t,x,u,\vartheta_{2})-f_{2}(t,x,u,\vartheta_{1})
+k(\vartheta_{2}-\vartheta_{1})\\[5pt]
&=f_{2}(t,x,u,\vartheta_{2})-f_{2}(t,x,u,\vartheta_{1})
+[k-d_{2}(1-\tau)](\vartheta_{2}-\vartheta_{1})\\[5pt]
&\geq 0,
\end{array}
\end{align*}
we can get $w(t,x)\geq 0$ in $\overline{\Omega}_{T_{0}}^{g,h}$
by the maximum principle for linear parabolic equation,
which implies $\Phi \vartheta_{2}\geq \Phi \vartheta_{1}$. Similarly, we can show that $\Phi \vartheta\leq \vartheta$ if $\vartheta$ is an upper solution, and $\Phi \vartheta\geq \vartheta$ if $\vartheta$ is a lower solution.

We then construct two sequences $\{v^{(n)}\}$ and $\{w^{(n)}\}$ by defining
$v^{(1)}=\Phi\bar{v},~v^{(n)}=\Phi v^{(n-1)},~w^{(1)}=\Phi\underline{v},~w^{(n)}=\Phi w^{(n-1)},
~n\geq 2$.
Thus,
$\underline{v}\leq w^{(1)}\leq w^{(2)}\leq \cdots
\leq w^{(n)}\leq v^{(n)}\leq \cdots\leq v^{(2)}
\leq v^{(1)}\leq \bar{v}$.
We conclude that the pointwise limits
\begin{align*}
w^{*}(t,x)=\lim_{n\rightarrow\infty}w^{(n)}(t,x),~
v^{*}(t,x)=\lim_{n\rightarrow\infty}v^{(n)}(t,x)
\end{align*}
exist at each point in $\Omega_{T_{0}}^{g,h}$ and
\begin{align*}
\underline{v}(t,x)\leq w^{*}(t,x)\leq v^{*}(t,x)\leq\bar{v}(t,x)
\quad \mbox{in}~\Omega_{T_{0}}^{g,h}.
\end{align*}

Similar as the proof of Theorem 2.4.6 in \cite{ylww11}, we can show that $v^{*}, w^{*}$ are classical solutions of $(\ref{2.5})$ and satisfy $v^{*}=w^{*}$. Moreover, the solution in $[\underline{v},\bar{v}]$
is unique.

\emph{Step 2.} It is easy to check that $\underline{v}=0$ and $\bar{v}=K_{2}$ are lower and upper solutions of $(\ref{2.5})$, respectively. Then there exists a unique solution $v$ satisfying $0<v\leq K_{2}$. Note that $f_{2}(t,x,u,v)$ satisfies the assumption (\textbf{f4}). Lemma $\ref{l2.2}$ implies that $v$ is unique solution of $(\ref{2.5})$.

We define
\begin{align*}
\begin{array}{rl}
\Omega:=
\Big\{(t,x):0<t\leq T_{0},~h(t)-M^{-1}<x<h(t)\Big\}
\end{array}
\end{align*}
and construct an auxiliary function
\begin{align*}
\psi(t,x)=K_{2}[2M(h(t)-x)-M^{2}(h(t)-x)^{2}].
\end{align*}
We will choose $M$ such that $\psi(t,x)\geq v(t,x)$ holds over $\Omega$.

Direct calculations show that, for $(t,x)\in \Omega$,
\begin{align*}
&\partial_{t}\psi=2K_{2}Mh^{\prime}(t)(1-M(h(t)-x))\geq 0,\\[5pt]
&-\partial_{xx}\psi=2K_{2}M^{2},~f_{2}(t,x,u,v)\leq \hat{L}v.
\end{align*}
It follows that
\begin{align*}
\begin{array}{rl}
&\partial_{t}\psi-d_{2}\left[\tau \partial_{xx}\psi
+(1-\tau)\left(\int_{g(t)}^{h(t)}J(x-y)\psi(t,y)dy-\psi\right)\right]\\[5pt]
&\geq 2d_{2}\tau K_{2}M^{2}-d_{2}(1-\tau)K_{2}\int_{g(t)}^{h(t)}J(x-y)dy\\[5pt]
&\geq 2d_{2}\tau K_{2}M^{2}-d_{2}(1-\tau)K_{2}\geq \hat{L}K_{2}\\[5pt]
&\geq \hat{L}v\geq \partial_{t}v-d_{2}\left[\tau \partial_{xx}v
+(1-\tau)\left(\int_{g(t)}^{h(t)}J(x-y)v(t,y)dy-v\right)\right]
\quad \mbox{in}~\Omega,
\end{array}
\end{align*}
if $M^{2}\geq \frac{\hat{L}+d_{2}(1-\tau)}{2d_{2}\tau}$. On the other hand,
\begin{align*}
\begin{array}{rl}
\psi(t,h(t)-M^{-1})=K_{2}\geq v(t,h(t)-M^{-1}),\quad
\psi(t,h(t))=0=v(t,h(t)).
\end{array}
\end{align*}
Choosing
\begin{align*}
\begin{array}{rl}
M:=\max\left\{\sqrt{\frac{\hat{L}+d_{2}(1-\tau)}{2d_{2}\tau}},
\frac{4\|v_{0}\|_{C^{1}([-h_{0},h_{0}])}}{3K_{2}}\right\},
\end{array}
\end{align*}
we can prove that
$v_{0}(x)\leq \psi(0,x)$ for $x\in [h_{0}-M^{-1},h_{0}]$.
Then we can apply Lemma $\ref{l2.2}$ to $\psi-v$ over $\Omega$ to deduce that
\begin{align*}
v(t,x)\leq \psi(t,x)\quad\mbox{for}~(t,x)\in \Omega.
\end{align*}
It then follows that
$v_{x}(t,h(t))\geq-2K_{2}M$.
Moreover, since $v(t,h(t))=0$ and $v>0$ in $\Omega_{T_{0}}^{g,h}$,
we have $v_{x}(t,h(t))<0$. The estimates for $v_{x}(t,g(t))$ can be similarly
obtained.
\end{proof}

Now, by approximation method we get the unique strong solution of $(\ref{2.5})$ provided that $g^{\prime}(t),h^{\prime}(t)$ and $u(t,x)$ are only continuous functions, which plays an important role in the proof of Lemma $\ref{l2.5}$ later.

\begin{lemma}\label{l2.4} Suppose that $(\ref{J})$ holds, $(g,h)\in \mathbb{G}_{T_{0}}^{h_{0}}\times\mathbb{H}_{T_{0}}^{h_{0}}$, $u\in C(\overline{\Omega}_{T_{0}}^{g,h})$, $f_{2}$ satisfies (\textbf{f1})-(\textbf{f4})
and $v_{0}$ satisfies $(\ref{1.2})$. Then the problem $(\ref{2.5})$ admits a unique solution $v\in W_{p}^{1,2}(\Omega_{T_{0}}^{g,h})\cap C^{\frac{1+\alpha}{2},1+\alpha}(\overline{\Omega}_{T_{0}}^{g,h})$ with any $p>3$. Moreover, $v$ satisfies $(\ref{2.6})$.
\end{lemma}

\begin{proof}
\emph{Step 1.} (Uniqueness)
Let
\begin{align*}
\begin{array}{rl}
\tilde{v}(t,z)=v(t,x(t,z)),\quad\tilde{f}(t,z,\tilde{u}, \tilde{v})=f(t,x(t,z),u(t,x(t,z)),v(t,x(t,z))),
\end{array}
\end{align*}
then the problem becomes
\begin{align*}\label{2.10}
\left\{\begin{array}{l}
\partial_{t}\tilde{v}=
d_{2}\tau\xi(t)\partial_{z}^{2}\tilde{v}+\eta(t,z)\partial_{z}\tilde{v}
+d_{2}(1-\tau)\left(\frac{h(t)-g(t)}{2}\int_{-1}^{1}J(\frac{h(t)-g(t)}{2}(z-s))
\tilde{v}(t,s)ds-\tilde{v}\right)\\[5pt]
\qquad\quad
+\tilde{f}_{2}(t,z,\tilde{u}, \tilde{v}), \quad (t,z)\in D_{T_{0}},\\[5pt]
\tilde{v}(t,-1)=\tilde{v}(t,1)=0,\quad t\in(0, T_{0}],\\[5pt]
\tilde{v}(0, z)=v_{0}(h_{0}z), \quad z\in[-1, 1].
\end{array}\right.
 \tag{2.10}
\end{align*}
Assume that $v_{i}(t,x)\in W_{p}^{1,2}(\Omega_{T_{0}}^{g,h})\cap C^{\frac{1+\alpha}{2},1+\alpha}(\overline{\Omega}_{T_{0}}^{g,h})$, $i=1,2$, are two solutions of $(\ref{2.5})$, then $\tilde{v}_{i}(t,z)=v_{i}(t,x(t,z))\in W_{p}^{1,2}(D_{T_{0}})\cap C^{\frac{1+\alpha}{2},1+\alpha}(\overline{D}_{T_{0}})$ are two solutions of $(\ref{2.10})$.
Let $\tilde{w}=\tilde{v}_{1}-\tilde{v}_{2}$, then $\tilde{w}$ satisfies
\begin{align*}\label{2.11}
\left\{\begin{array}{l}
\partial_{t}\tilde{w}=
d_{2}\tau\xi(t)\partial_{z}^{2}\tilde{w}+\eta(t,z)\partial_{z}\tilde{w}
+d_{2}(1-\tau)\left(\frac{h(t)-g(t)}{2}\int_{-1}^{1}J(\frac{h(t)-g(t)}{2}(z-s))
\tilde{w}(t,s)ds-\tilde{w}\right)\\[5pt]
\qquad\quad
+\tilde{f}_{2}(t,z,\tilde{u}, \tilde{v}_{1})
-\tilde{f}_{2}(t,z,\tilde{u}, \tilde{v}_{2}), \quad (t,z)\in D_{T_{0}},\\[5pt]
\tilde{w}(t,-1)=\tilde{w}(t,1)=0,\quad t\in(0, T_{0}],\\[5pt]
\tilde{w}(0, z)=0, \quad z\in[-1, 1].
\end{array}\right.
 \tag{2.11}
\end{align*}
Multiplying the equation in $(\ref{2.11})$ by $\tilde{w}\chi_{[0,t]}$, where
$\chi_{[0,t]}$ is the characteristic function in $[0,t]$ with any $0<t\leq T_{0}$, and then integrating over $(0,T_{0}]\times [-1,1]$ gives
\begin{align*}
\begin{array}{rl}
&\frac{1}{2}\int_{-1}^{1}\tilde{w}^{2}(t,z)\Big|_{0}^{t}dz\\[5pt]
&=-d_{2}\tau\int_{0}^{t}\int_{-1}^{1}\xi(t)(\partial_{z}\tilde{w})^{2}dzdt
+\int_{0}^{t}\int_{-1}^{1}\eta(t,z)\tilde{w}\partial_{z}\tilde{w}dzdt\\[5pt]
&\quad+d_{2}(1-\tau)\int_{0}^{t}\int_{-1}^{1}\left(\frac{h(t)-g(t)}{2}\int_{-1}^{1}J(\frac{h(t)-g(t)}{2}(z-s))
\tilde{w}(t,s)ds-\tilde{w}\right)\tilde{w}dzdt\\[5pt]
&\quad
+\int_{0}^{t}\int_{-1}^{1}[\tilde{f}_{2}(t,z,\tilde{u},\tilde{v}_{1})
-\tilde{f}_{2}(t,z,\tilde{u},\tilde{v}_{2})]\tilde{w}dzdt.
\end{array}
\end{align*}
By the Young's inequality with $0<\varepsilon<\frac{4d_{2}\tau}{(h(T_{0})-g(T_{0}))^{2}}$,
\begin{align*}
\begin{array}{rl}
\int_{0}^{t}\int_{-1}^{1}\eta(t,z)\tilde{w}\partial_{z}\tilde{w}dzdt
\leq \varepsilon \int_{0}^{t}\int_{-1}^{1}(\partial_{z}\tilde{w})^{2}dzdt
+C(\varepsilon)\int_{0}^{t}\int_{-1}^{1}\tilde{w}^{2}dzdt.
\end{array}
\end{align*}
By the continuity of $J$ and H\"{o}lder inequality,
\begin{align*}
\begin{array}{rl}
&d_{2}(1-\tau)\int_{0}^{t}\int_{-1}^{1}\left(\frac{h(t)-g(t)}{2}\int_{-1}^{1}J(\frac{h(t)-g(t)}{2}(z-s))
\tilde{w}(t,s)ds-\tilde{w}(t,z)\right)\tilde{w}(t,z)dzdt\\[5pt]
&\leq d_{2}(1-\tau)C\int_{0}^{t}(\int_{-1}^{1}|\tilde{w}(t,z)|dz)^{2}dt
-d_{2}(1-\tau)\int_{0}^{t}\int_{-1}^{1}\tilde{w}^{2}dzdt\\[5pt]
&\leq d_{2}(1-\tau)C_{1}\int_{0}^{t}\int_{-1}^{1}\tilde{w}^{2}dzdt.
\end{array}
\end{align*}
By the Lipschitz continuity of $f_{2}$ with respect to $\tilde{v}$,
\begin{align*}
\begin{array}{rl}
\int_{0}^{t}\int_{-1}^{1}[\tilde{f}_{2}(t,z,\tilde{u},\tilde{v}_{1})
-\tilde{f}_{2}(t,z,\tilde{u},\tilde{v}_{2})]\tilde{w}(t,z)dzdt
\leq L\int_{0}^{t}\int_{-1}^{1}\tilde{w}^{2}dzdt.
\end{array}
\end{align*}
Combining the above estimates, we have
\begin{align*}
\begin{array}{rl}
\int_{-1}^{1}\tilde{w}^{2}(t,z)dz
\leq C\int_{0}^{t}\int_{-1}^{1}\tilde{w}^{2}dzdt.
\end{array}
\end{align*}
By the Gronwall's inequality, we know $\int_{0}^{t}\int_{-1}^{1}\tilde{w}^{2}dzdt=0$, which implies that $\tilde{w}=0$, a.e. in $(0,t]\times [-1,1]$.
Since $t\in(0,T_{0}]$ is arbitrary and $\tilde{w}\in C(\overline{D}_{T_{0}})$, we can obtain $\tilde{w}=0$ for all $(t,z)$ in $[0,T_{0}]\times [-1,1]$, which implies the uniqueness of solution.

\emph{Step 2.} (Existence) For any $(g,h)\in \mathbb{G}_{T_{0}}^{h_{0}}\times\mathbb{H}_{T_{0}}^{h_{0}}$, we can find some sequences $(g_{n},h_{n})\in \widehat{\mathbb{G}}_{T_{0}}^{h_{0}}\times\widehat{\mathbb{H}}_{T_{0}}^{h_{0}}$ such that $g_{n}\rightarrow g$ and $h_{n}\rightarrow h$ in $C^{1}([0,T_{0}])$. Moreover, for every $u(t,x)\in C(\overline{\Omega}_{T_{0}}^{g,h})$, we can obtain $\tilde{u}(t,z)=u(t,x(t,z))\in C(\overline{D}_{T_{0}})$ and find some sequence $\tilde{u}_{n}\in C^{\frac{\alpha}{2},\alpha}(\overline{D}_{T_{0}})$ such that
$\tilde{u}_{n}\rightarrow \tilde{u}$ in $C(\overline{D}_{T_{0}})$. Taking
$u_{n}(t,x)=\tilde{u}_{n}(t,\frac{2x-g_{n}(t)-h_{n}(t)}{h_{n}(t)-g_{n}(t)})$, we know
$u_{n}\in C^{\frac{\alpha}{2},\alpha}(\overline{\Omega}_{T_{0}}^{g_{n},h_{n}})$.

Consider the approximate problem
\begin{align*}\label{2.12}
\left\{\begin{array}{l}
\partial_{t}v=
d_{2}\left[\tau \partial_{x}^{2}v
+(1-\tau)\left(\int_{g_{n}(t)}^{h_{n}(t)}J(x-y)v(t,y)dy-v\right)\right]\\[5pt]
\qquad\quad+f_{2}(t,x,u_{n},v), \quad (t,x)\in\Omega_{T_{0}}^{g_{n},h_{n}},\\[5pt]
v(t,g_{n}(t))=v(t,h_{n}(t))=0,\quad t\in(0, T_{0}],\\[5pt]
v(0, x)=v_{0}(x), \quad x\in[-h_{0}, h_{0}].
\end{array}\right.
 \tag{2.12}
\end{align*}
By Lemma $\ref{l2.3}$, we know $(\ref{2.12})$ has a unique classical solution $v_{n}\in C^{1+\frac{\alpha}{2},2+\alpha}(\Omega_{T_{0}}^{g_{n},h_{n}})$, and satisfies
\begin{align*}
&0<v_{n}\leq K_{2}\quad\mbox{for}~(t,x)\in\Omega_{T_{0}}^{g_{n},h_{n}},\\[3pt]
&0<-\partial_{x}v_{n}(t,h_{n}(t)), \partial_{x}v_{n}(t,g_{n}(t))\leq K_{3}
\quad\mbox{for}~t\in(0,T_{0}].
\end{align*}

Let $\tilde{v}_{n}(t,z)=v_{n}(t,x_{n}(t,z))$ and
\begin{align*}
\begin{array}{rl}
\tilde{f}_{2}(t,z,\tilde{u}_{n}, \tilde{v}_{n})
=f_{2}(t,x_{n}(t,z),u_{n}(t,x_{n}(t,z)),v_{n}(t,x_{n}(t,z)))
\end{array}
\end{align*}
with
\begin{align*}
\begin{array}{rl}
x_{n}(t,z)=\frac{(h_{n}(t)-g_{n}(t))z+h_{n}(t)+g_{n}(t)}{2},
\end{array}
\end{align*}
then $\tilde{v}_{n}(t,z)\in C^{1+\frac{\alpha}{2},2+\alpha}(D_{T_{0}})$ is the unique solution of
\begin{align*}\label{2.13}
\left\{\begin{array}{l}
\partial_{t}\tilde{v}_{n}=
d_{2}\tau\xi_{n}(t)\partial_{z}^{2}\tilde{v}_{n}+\eta_{n}(t,z)\partial_{z}\tilde{v}_{n}\\[5pt]
\qquad\quad
+d_{2}(1-\tau)\left(\frac{h_{n}(t)-g_{n}(t)}{2}\int_{-1}^{1}J(\frac{h_{n}(t)-g_{n}(t)}{2}(z-s))
\tilde{v}_{n}(t,s)ds-\tilde{v}_{n}\right)\\[5pt]
\qquad\quad
+\tilde{f}_{2}(t,z,\tilde{u}_{n}, \tilde{v}_{n}), \quad (t,z)\in D_{T_{0}},\\[5pt]
\tilde{v}_{n}(t,-1)=\tilde{v}_{n}(t,1)=0,\quad t\in(0, T_{0}],\\[5pt]
\tilde{v}_{n}(0, z)=v_{0}(h_{0}z), \quad z\in[-1, 1],
\end{array}\right.
 \tag{2.13}
\end{align*}
and satisfies
\begin{align*}\label{2.14}
\begin{array}{rl}
&0<\tilde{v}_{n}\leq K_{2}\quad\mbox{in}~D_{T_{0}},\\[5pt]
&0<-\frac{2}{h_{n}(t)-g_{n}(t)}\partial_{z}\tilde{v}_{n}(t,1),~
\frac{2}{h_{n}(t)-g_{n}(t)}\partial_{z}\tilde{v}_{n}(t,-1)\leq K_{3}
\quad\mbox{for}~t\in(0,T_{0}].
\end{array}
 \tag{2.14}
\end{align*}

Let
\begin{align*}
\begin{array}{rl}
g(t,z):=
d_{2}(1-\tau)\left(\frac{h_{n}(t)-g_{n}(t)}{2}\int_{-1}^{1}J(\frac{h_{n}(t)-g_{n}(t)}{2}(z-s))
\tilde{v}_{n}(t,s)ds\right)
+\tilde{f}_{2}(t,z,\tilde{u}_{n}, \tilde{v}_{n}),
\end{array}
\end{align*}
we know $g\in L^{\infty}(D_{T_{0}})$. Applying the $L^{p}$ theory for linear parabolic equations
to $(\ref{2.13})$, we have the solution $\tilde{v}_{n}$ satisfies
$\|\tilde{v}_{n}\|_{W_{p}^{1,2}(D_{T_{0}})}\leq C$,
where $C$ is independent of $n$. By the weak compactness of the bounded set in $W_{p}^{1,2}(D_{T_{0}})$ and $\mathring{W}_{p}^{1,1}(D_{T_{0}})$ and the compactly imbedding theorem ($W_{p}^{1,1}(D_{T_{0}})\hookrightarrow\hookrightarrow
L^{p}(D_{T_{0}})$), there exists a subsequence, still denoted
by $\{\tilde{v}_{n}\}$, such that $\tilde{v}_{n}\rightharpoonup \tilde{v}$ in $W_{p}^{1,2}(D_{T_{0}})\cap\mathring{W}_{p}^{1,1}(D_{T_{0}})$, $\partial_{z}\tilde{v}_{n}\rightarrow \partial_{z}\tilde{v}$ in $L^{p}(D_{T_{0}})$ and $\tilde{v}_{n}\rightarrow \tilde{v}$ in $L^{p}(D_{T_{0}})$, which implies that $\tilde{v}\in W_{p}^{1,2}(D_{T_{0}})\cap
\mathring{W}_{p}^{1,1}(D_{T_{0}})$ is the strong solution of $(\ref{2.10})$. By the Sobolev imbedding theorem, $\tilde{v}\in C^{\frac{1+\alpha}{2},1+\alpha}(\overline{D}_{T_{0}})$.

Note that $\tilde{v}_{n}$ satisfies $(\ref{2.14})$. From the fact
$\partial_{z}\tilde{v}_{n}\rightarrow \partial_{z}\tilde{v}$, $\tilde{v}_{n}\rightarrow \tilde{v}$ in $L^{p}(D_{T_{0}})$ (then a.e. in $D_{T_{0}}$) and $\tilde{v}\in C^{\frac{1+\alpha}{2},1+\alpha}(\overline{D}_{T_{0}})$,
we have $0<\tilde{v}\leq K_{2}$ in $D_{T_{0}}$ and
$0<-\frac{2}{h(t)-g(t)}\partial_{z}\tilde{v}(t,1),
\frac{2}{h(t)-g(t)}\partial_{z}\tilde{v}(t,-1)\leq K_{3}$
for $t\in(0,T_{0}]$. Then $v(t,x)=\tilde{v}(t,z(t,x))$ satisfies
$(\ref{2.6})$, which completes the proof.
\end{proof}

In the following lemma, we prove the well-posedness for $(\ref{2.1})$ with any fixed $(g,h)\in \mathbb{G}_{T_{0}}^{h_{0}}\times \mathbb{H}_{T_{0}}^{h_{0}}$ by the fixed point theorem. Denote
\begin{align*}
\begin{array}{rl}
&\mathbb{X}_{T_{0}}^{1}
:=\Big\{u\in C(\overline{\Omega}_{T_{0}}^{g,h}):~0\leq u\leq K_{1}, u(0,x)=u_{0}(x), u(t,g(t))=u(t,h(t))=0\Big\}, \\[5pt]
&\mathbb{X}_{T_{0}}^{2}
:=\Big\{v\in C(\overline{\Omega}_{T_{0}}^{g,h}):~0\leq v\leq K_{2}, v(0,x)=v_{0}(x), v(t,g(t))=v(t,h(t))=0\Big\}, \\[5pt]
&\mathbb{X}_{T_{0}}^{g,h}:=\mathbb{X}_{T_{0}}^{1}\times\mathbb{X}_{T_{0}}^{2}.
\end{array}
\end{align*}

\begin{lemma}\label{l2.5}  For any $T_{0}>0$ and $(g,h)\in \mathbb{G}_{T_{0}}^{h_{0}}\times \mathbb{H}_{T_{0}}^{h_{0}}$, the problem
\begin{align*}\label{2.15}
\left\{\begin{array}{l}
\partial_{t}u=d_{1}\left(\int_{g(t)}^{h(t)}J(x-y)u(t,y)dy-u\right)+f_{1}(t,x,u,v),~(t,x)\in \Omega_{T_{0}}^{g,h},\\[5pt]
\partial_{t}v=d_{2}\Big[\tau \partial_{x}^{2}v+(1-\tau)\left(\int_{g(t)}^{h(t)}J(x-y)v(t,y)dy-v\right)\Big]\\[5pt]
\qquad\quad+f_{2}(t,x,u,v),\quad (t,x)\in \Omega_{T_{0}}^{g,h},\\[5pt]
u(t, g(t))=u(t,h(t))=v(t, g(t))=v(t,h(t))=0, \quad t\in [0,T_{0}],\\[5pt]
u(0,x)=u_0(x), v(0, x)=v_0(x), \quad x\in[-h_{0},h_{0}]
\end{array}\right.
\tag{2.15}
\end{align*}
admits a unique solution $(u, v)\in \mathbb{X}_{T_{0}}^{g,h}$,
and $(u, v)$ satisfy
\begin{align*}\label{2.16}
\begin{array}{rl}
&0<u\leq K_{1}, 0<v\leq K_{2}\quad \mbox{in}~\Omega_{T_{0}}^{g,h},\\[5pt]
&0<-v_{x}(t,h(t)), v_{x}(t,g(t))\leq K_{3}\quad \mbox{in}~(0,T_{0}].
\end{array}
\tag{2.16}
\end{align*}
Moreover, $v\in W_{p}^{1,2}(\Omega_{T_{0}}^{g,h})\cap C^{\frac{1+\alpha}{2},1+\alpha}(\overline{\Omega}_{T_{0}}^{g,h})$
with any $p>3$.
\end{lemma}

\begin{proof} For $u^{*}\in \mathbb{X}_{s}^{1}$ with $0<s\leq T_{0}$, from Lemma $\ref{l2.4}$ we know that the initial-boundary value problem $(\ref{2.5})$ with $(u,T_{0})$ replaced by $(u^{*},s)$ admits a unique solution $v\in\mathbb{X}_{s}^{2}$. For such $v\in\mathbb{X}_{s}^{2}$, we consider
\begin{align*}
\left\{\begin{array}{l}
\partial_{t}u=d_{1}\left(\int_{g(t)}^{h(t)}J(x-y)u(t,y)dy-u\right)+f_{1}(t,x,u,v), \quad (t,x)\in \Omega_{T_{0}}^{g,h},\\[5pt]
u(t, g(t))=u(t,h(t))=0, \quad t\in [0,T_{0}],\\[5pt]
u(0,x)=u_0(x), \quad x\in[-h_{0},h_{0}].
\end{array}\right.
\end{align*}
By Lemma $2.3$ in \cite{cdll19}, it admits a unique solution $u\in \mathbb{X}_{s}^{1}$.
We define a mapping
$\mathcal{F}_{s}: \mathbb{X}_{s}^{1}\rightarrow\mathbb{X}_{s}^{1}$
by $\mathcal{F}_{s}u^{*}=u$. If $\mathcal{F}_{s}u^{*}=u^{*}$, then
$(u^{*},v)$ solves $(\ref{2.15})$ with $T_{0}$ replaced by $s$.

Next, we shall prove that
$\mathcal{F}_{s}$ has a fixed point in $\mathbb{X}_{s}^{1}$ provided that $s$ is small
enough. For $i=1,2$, we assume $u_{i}^{*}\in \mathbb{X}_{s}^{1}$, $u_{i}=\mathcal{F}_{s}u_{i}^{*}$, and $v_{i}$ be the unique solution of $(\ref{2.5})$ with $(u,T_{0})$ replaced by $(u_{i}^{*},s)$.
Denote $\theta^{*}=u_{1}^{*}-u_{2}^{*}$, $\theta=u_{1}-u_{2}$ and $w=v_{1}-v_{2}$.
Note that $w$ satisfies
\begin{align*}
\left\{\begin{array}{l}
\partial_{t}w=d_{2}\left[\tau \partial_{x}^{2}w+(1-\tau)\left(\int_{g(t)}^{h(t)}J(x-y)w(t,y)dy-w\right)\right]\\[5pt]
\qquad\quad+a_{0}(t,x)w+b_{0}(t,x)\theta^{*}, \quad (t,x)\in \Omega_{s}^{g,h},\\[5pt]
w(t, g(t))=w(t,h(t))=0, \quad t\in [0,s],\\[5pt]
w(0, x)=0, \quad x\in[-h_{0},h_{0}],
\end{array}\right.
\end{align*}
where
\begin{align*}
\begin{array}{rl}
a_{0}(t,x)=\int_{0}^{1}f_{2,v}(t,x,u_{1}^{*},v_{2}+(v_{1}-v_{2})\tau)d\tau,\\[3pt]
b_{0}(t,x)=\int_{0}^{1}f_{2,u}(t,x,u_{2}^{*}+(u_{1}^{*}-u_{2}^{*})\tau,v_{2})d\tau.
\end{array}
\end{align*}
Let
$\tilde{\theta}^{*}(t,z)=\theta^{*}(t,x(t,z)),
\tilde{w}(t,z)=w(t,x(t,z)),
\tilde{a}_{0}(t,z)=a_{0}(t,x(t,z)),
\tilde{b}_{0}(t,z)=b_{0}(t,x(t,z))$.
It is easy to see that $\tilde{w}$ satisfies
\begin{align*}
\left\{\begin{array}{l}
\partial_{t}\tilde{w}=d_{2}\tau\xi(t)\partial_{z}^{2}\tilde{w}
+\eta(t,z)\partial_{z}\tilde{w}+[\tilde{a}_{0}(t,z)-d_{2}(1-\tau)]\tilde{w}\\[5pt]
\qquad
+d_{2}(1-\tau)\frac{h(t)-g(t)}{2}\int_{-1}^{1}J(\frac{h(t)-g(t)}{2}(z-s))
\tilde{w}(t,s)ds
+\tilde{b}_{0}(t,z)\tilde{\theta}^{*}, \quad (t,z)\in D_{s},\\[5pt]
\tilde{w}(t, -1)=\tilde{w}(t,1)=0, \quad t\in [0,s],\\[5pt]
\tilde{w}(0, z)=0, \quad z\in[-1,1].
\end{array}\right.
\end{align*}
By the $L^{p}$ theory for linear parabolic equation, we have
\begin{align*}
\begin{array}{rl}
\|\tilde{w}\|_{W_{p}^{1,2}(D_{s})}
&\leq C\left(\left\|\frac{h(t)-g(t)}{2}\int_{-1}^{1}J(\frac{h(t)-g(t)}{2}(z-s))
\tilde{w}(t,s)ds\right\|_{L^{p}(D_{s})}
+\|\tilde{\theta}^{*}\|_{L^{p}(D_{s})}\right)\\[3pt]
&\leq C(\|\tilde{w}\|_{C(\overline{D}_{s})}
\left\|\int_{\frac{h(t)-g(t)}{2}(z-1)}^{\frac{h(t)-g(t)}{2}(z+1)}J(y)dy\right\|_{L^{p}(D_{s})}
+\|\tilde{\theta}^{*}\|_{L^{p}(D_{s})})\\[3pt]
&\leq C(\|\tilde{w}\|_{C(\overline{D}_{s})}(2s)^{\frac{1}{p}}
+\|\tilde{\theta}^{*}\|_{L^{p}(D_{s})}).
\end{array}
\end{align*}
From the proof of Theorem 1.1 in \cite{w19}, we know the H\"{o}lder
semi-norm $[\tilde{w}]_{C^{\frac{\alpha}{2},\alpha}(\overline{D}_{s})}
\leq C^{\prime}\|\tilde{w}\|_{W_{p}^{1,2}(D_{s})}$, where $C^{\prime}$ is independent
of $\frac{1}{s}$. Thus,
\begin{align*}
|\tilde{w}(t,z)|
=|\tilde{w}(t,z)-\tilde{w}(0,z)|
\leq [\tilde{w}]_{C^{\frac{\alpha}{2},\alpha}(\overline{D}_{s})}t^{\frac{\alpha}{2}}
\leq C^{\prime}\|\tilde{w}\|_{W_{p}^{1,2}(D_{s})}t^{\frac{\alpha}{2}},
\end{align*}
which implies that
\begin{align*}
\|\tilde{w}\|_{C(\overline{D}_{s})}
\leq C^{\prime}\|\tilde{w}\|_{W_{p}^{1,2}(D_{s})}s^{\frac{\alpha}{2}}.
\end{align*}
Choosing $s$ small such that $CC^{\prime}(2s)^{\frac{1}{p}}s^{\frac{\alpha}{2}}<\frac{1}{2}$, we have
\begin{align*}
\begin{array}{rl}
\|\tilde{w}\|_{W_{p}^{1,2}(D_{s})}
\leq 2C\|\tilde{\theta}^{*}\|_{L^{p}(D_{s})}
\leq 2C(2s)^{\frac{1}{p}}\|\tilde{\theta}^{*}\|_{C(\overline{D}_{s})}
=2C(2s)^{\frac{1}{p}}\|\theta^{*}\|_{C(\overline{\Omega}_{s}^{g,h})}.
\end{array}
\end{align*}
Similar to the proof of Lemma 2.3 (Step 3) in \cite{waw182}, we can choose $s$ small enough
such that
\begin{align*}
\begin{array}{rl}
\|\theta\|_{C(\overline{\Omega}_{s}^{g,h})}
\leq\frac{1}{2}\|\theta^{*}\|_{C(\overline{\Omega}_{s}^{g,h})}.
\end{array}
\end{align*}
By the contraction mapping theorem, we know that $\mathcal{F}_{s}$
has a unique fixed point $u\in \mathbb{X}_{s}^{1}$.

Following the arguments in the proof of Lemma 2.3 (Step 5) in \cite{waw182}, we can show that the unique solution $(u,v)$ of $(\ref{2.15})$ can be extended to $\Omega_{T_{0}}^{g,h}$ and $(u,v)\in \mathbb{X}_{T_{0}}^{g,h}$. The estimates of $v_{x}(t,h(t)), v_{x}(t,g(t))$ and the regularity of $v$
have been established in Lemma $\ref{l2.4}$.
\end{proof}

\noindent\textbf{Proof of Theorem $\ref{t2.1}$.} By Lemma $\ref{l2.5}$, for any $T_{0}>0$ and
$(g,h)\in \mathbb{G}_{T_{0}}^{h_{0}}\times \mathbb{H}_{T_{0}}^{h_{0}}$,
we can find a unique $(u,v)\in \mathbb{X}_{T_{0}}^{g,h}$ that solves $(\ref{2.15})$, and $(\ref{2.16})$ holds. For $0<t\leq T_{0}$, define the mapping
\begin{align*}
\begin{array}{rl}
\mathcal{G}(g,h)=(\tilde{g},\tilde{h})
\end{array}
\end{align*}
by
\begin{align*}
\begin{array}{rl}
\tilde{h}(t)
&=h_{0}-\mu \int_{0}^{t}v_{x}(\tau, h(\tau))d\tau
+\rho_{1}\int_{0}^{t}\int_{g(\tau)}^{h(\tau)}\int_{h(\tau)}^{\infty}J(x-y)u(\tau,x)dydxd\tau\\[3pt]
&\quad+\rho_{2}\int_{0}^{t}\int_{g(\tau)}^{h(\tau)}\int_{h(\tau)}^{\infty}J(x-y)v(\tau,x)dydxd\tau,\\[3pt]
\tilde{g}(t)
&=-h_{0}-\mu \int_{0}^{t}v_{x}(\tau, g(\tau))d\tau
-\rho_{1}\int_{0}^{t}\int_{g(\tau)}^{h(\tau)}\int_{-\infty}^{g(\tau)}J(x-y)u(\tau,x)dydxd\tau\\[3pt]
&\quad-\rho_{2}\int_{0}^{t}\int_{g(\tau)}^{h(\tau)}\int_{-\infty}^{g(\tau)}J(x-y)v(\tau,x)dydxd\tau.
\end{array}
\end{align*}
To prove this theorem, we will show that if $T_{0}$ is sufficiently small, then $\mathcal{G}$ maps a suitable closed subset $\Sigma_{T_{0}}$ of $\mathbb{G}_{T_{0}}^{h_{0}}\times \mathbb{H}_{T_{0}}^{h_{0}}$ into itself and is a contraction mapping.
The proof can be obtained by using similar arguments as that of Theorem 2.1 in \cite{dwz19,waw182},
here we omit the details. \hfill $\Box$

\section{Comparison principle and some eigenvalue problems}

In this section, we first give two comparison principles for $(\ref{1.1})$ and $T$-periodic nonlocal evolution equation $(\ref{3.2})$, and then investigate the existence and properties of principle eigenvalue of some eigenvalue problems. These results will play an important role in later sections.

\subsection{The comparison principle}

In this subsection, we discuss the comparison principle for $(\ref{1.1})$.

\begin{lemma}\label{l3.1}
(The Comparison Principle) Suppose that $T_{0}\in(0, \infty)$, $\bar{g}, \bar{h}\in C^1([0, T_{0}])$, $\bar{u}\in C(\overline{\Omega}_{T_{0}}^{\bar{g}, \bar{h}})$, $\bar{v}\in C^{1,2}(\Omega_{T_{0}}^{\bar{g},\bar{h}})\cap C(\overline{\Omega}_{T_{0}}^{\bar{g},\bar{h}})$,
and $(\bar{u},\bar{v},\bar{g},\bar{h})$ satisfy
\begin{align*}\label{3.1}
\left\{\begin{array}{l}
\partial_{t}\bar{u}\geq d_{1}\left(\int_{\bar{g}(t)}^{\bar{h}(t)}J(x-y)\bar{u}(t,y)dy-\bar{u}\right)
+\bar{u}(a(t)-\bar{u}), \quad (t,x)\in \Omega_{T_{0}}^{\bar{g},\bar{h}},\\[5pt]
\partial_{t}\bar{v}\geq
d_{2}\left[\tau \partial_{x}^{2}\bar{v}+(1-\tau)\left(\int_{\bar{g}(t)}^{\bar{h}(t)}J(x-y)\bar{v}(t,y)dy-\bar{v}\right)\right]
+\bar{v}(c(t)-\bar{v}), \quad (t,x)\in \Omega_{T_{0}}^{\bar{g},\bar{h}},\\[5pt]
\bar{u}(t, \bar{g}(t))\geq0, \bar{u}(t,\bar{h}(t))\geq0, \quad 0<t\leq T_{0},\\[5pt]
\bar{v}(t, \bar{g}(t))=0, \bar{v}(t,\bar{h}(t))=0, \quad 0<t\leq T_{0},\\[5pt]
\bar{h}'(t)\geq-\mu \bar{v}_{x}(t, \bar{h}(t))
+\rho_{1}\int_{\bar{g}(t)}^{\bar{h}(t)}\int_{\bar{h}(t)}^{\infty}J(x-y)\bar{u}(t,x)dydx\\[5pt]
\qquad\quad+\rho_{2}\int_{\bar{g}(t)}^{\bar{h}(t)}\int_{\bar{h}(t)}^{\infty}J(x-y)\bar{v}(t,x)dydx, \quad 0<t\leq T_{0},\\[5pt]
\bar{g}'(t)\leq-\mu \bar{v}_{x}(t, \bar{g}(t))
-\rho_{1}\int_{\bar{g}(t)}^{\bar{h}(t)}\int_{-\infty}^{\bar{g}(t)}J(x-y)\bar{u}(t,x)dydx\\[5pt]
\qquad\quad-\rho_{2}\int_{\bar{g}(t)}^{\bar{h}(t)}\int_{-\infty}^{\bar{g}(t)}J(x-y)\bar{v}(t,x)dydx, \quad 0<t\leq T_{0},\\[5pt]
\bar{u}(0,x)\geq u_0(x), \bar{v}(0, x)\geq v_0(x), \quad |x|\leq h_{0},\\[5pt]
\bar{h}(0)\geq h_{0},~\bar{g}(0)\leq-h_{0}.
\end{array}\right.
\tag{3.1}
\end{align*}
Let $(u,v,g,h)$ be the unique solution of $(\ref{1.1})$, then
\begin{align*}
\begin{array}{l}
g(t)\geq \bar{g}(t),~
h(t)\leq \bar{h}(t)~ \mbox{in}~(0,T_{0}],\\
u(t,x)\leq \bar{u}(t,x),~
v(t,x)\leq \bar{v}(t,x)~
\mbox{for}~(t,x)\in \overline{\Omega}_{T_{0}}^{g,h}.
\end{array}
\end{align*}
\end{lemma}

\begin{proof}
Thanks to Lemma 2.2 in \cite{cdll19} and Lemma $\ref{l2.2}$, one sees that $\bar{u},\bar{v}>0$
for $(t,x)\in \Omega_{T_{0}}^{\bar{g},\bar{h}}$.

We first consider the case $\bar{h}(0)>h_{0},~\bar{g}(0)<-h_{0}$. Then $\bar{h}(t)>h(t),~\bar{g}(t)<g(t)$
hold true for small $t>0$. We claim that $\bar{h}(t)>h(t),~\bar{g}(t)<g(t)$ for all $t\in (0,T_{0}]$.
In fact, if this is not true, there exists $t_{1}\leq T_{0}$ such that
\begin{align*}
\bar{h}(t)>h(t),~\bar{g}(t)<g(t) \quad \mbox{for}~t\in(0,t_{1})\quad
\mbox{and}\quad [\bar{h}(t_{1})-h(t_{1})][\bar{g}(t_{1})-g(t_{1})]=0.
\end{align*}
Without loss of generality, we may assume that
\begin{align*}
\bar{g}(t_{1})\leq g(t_{1})\quad
\mbox{and}\quad
\bar{h}(t_{1})=h(t_{1}).
\end{align*}
Thus, $\bar{h}^{\prime}(t_{1})\leq h^{\prime}(t_{1})$.
Since $\bar{v}(0,x)\geq v_0(x)$ for $x\in[-h_{0},h_{0}]$,
$\bar{v}(t,g(t))\geq0=v(t,g(t))$ and $\bar{v}(t,h(t))\geq0=v(t,h(t))$ for $t\in (0,t_{1}]$,
by applying Lemma $\ref{l2.2}$, we have $\bar{v}>v$ in $\Omega_{t_{1}}^{\bar{g},\bar{h}}$.
Moreover, by the fact that $\bar{v}(t_{1},h(t_{1}))=\bar{v}(t_{1},\bar{h}(t_{1}))=0=v(t_{1},h(t_{1}))$,
we deduce that $\bar{v}_{x}(t_{1},h(t_{1}))<v_{x}(t_{1},h(t_{1}))$. Similarly, using Lemma 2.2 in \cite{cdll19}, we can obtain $\bar{u}>u$ in $\Omega_{t_{1}}^{\bar{g},\bar{h}}$.
It follows that
\begin{align*}
\begin{array}{rl}
\bar{h}^{\prime}(t_{1})
&\geq -\mu \bar{v}_{x}(t_{1}, \bar{h}(t_{1}))
+\rho_{1}\int_{\bar{g}(t_{1})}^{\bar{h}(t_{1})}\int_{\bar{h}(t_{1})}^{\infty}J(x-y)\bar{u}(t_{1},x)dydx\\[3pt]
&\quad+\rho_{2}\int_{\bar{g}(t_{1})}^{\bar{h}(t_{1})}\int_{\bar{h}(t_{1})}^{\infty}J(x-y)\bar{v}(t_{1},x)dydx\\[3pt]
&\geq -\mu \bar{v}_{x}(t_{1}, h(t_{1}))
+\rho_{1}\int_{g(t_{1})}^{h(t_{1})}\int_{h(t_{1})}^{\infty}J(x-y)\bar{u}(t_{1},x)dydx\\[3pt]
&\quad+\rho_{2}\int_{g(t_{1})}^{h(t_{1})}\int_{h(t_{1})}^{\infty}J(x-y)\bar{v}(t_{1},x)dydx\\[3pt]
&>-\mu v_{x}(t_{1}, h(t_{1}))
+\rho_{1}\int_{g(t_{1})}^{h(t_{1})}\int_{h(t_{1})}^{\infty}J(x-y)u(t_{1},x)dydx\\[3pt]
&\quad+\rho_{2}\int_{g(t_{1})}^{h(t_{1})}\int_{h(t_{1})}^{\infty}J(x-y)v(t_{1},x)dydx\\[3pt]
&=h^{\prime}(t_{1}),
\end{array}
\end{align*}
which is a contradiction. Hence, $h(t)<\bar{h}(t)$, $g(t)>\bar{g}(t)$ for all $t\in (0,T_{0}]$,
and $\bar{u}(t,x)>u(t,x)$, $\bar{v}(t,x)>v(t,x)$ in $\Omega_{T_{0}}^{g,h}$.

For the general case that $\bar{h}(0)\geq h_{0},~\bar{g}(0)\leq-h_{0}$, we can adopt the same method as the proof Lemma 5.1 in \cite{gw12}.
\end{proof}

\begin{remark}\label{r3.2} From the proof of Lemma $\ref{l3.1}$, we can see that the conditions $\bar{v}(t, \bar{g}(t))=0$, $\bar{v}(t,\bar{h}(t))=0$ are necessary in deriving the contradiction from the relationship
between $\bar{h}^{\prime}(t)$ and $h^{\prime}(t)$. If $\tau=0$, as considered in \cite{dwz19}, then the expressions of $h^{\prime}(t), g^{\prime}(t)$ in (1.1) and $\bar{h}^{\prime}(t)$, $\bar{g}^{\prime}(t)$ in $(\ref{3.1})$ do not include the terms $-\mu v_{x}(t,h(t))$, $-\mu v_{x}(t,g(t))$ and $-\mu \bar{v}_{x}(t,\bar{h}(t))$, $-\mu \bar{v}_{x}(t,\bar{g}(t))$, respectively, in such case the conditions $\bar{v}(t, \bar{g}(t))=0$, $\bar{v}(t,\bar{h}(t))=0$ can be
weaken into $\bar{v}(t, \bar{g}(t))\geq0$, $\bar{v}(t,\bar{h}(t))\geq0$.
\end{remark}

In the following, we establish a comparison principle for the following $T$-periodic nonlocal evolution equation
\begin{align*}\label{3.2}
\left\{\begin{array}{l}
u_{t}=d_{1}[\int_{\Omega}J(x-y)u(t,y)dy-u(t,x)]+u(a(t)-u),
\quad (t,x)\in \mathbb{R}\times \overline{\Omega},\\[3pt]
u(0,x)=u(T,x), \quad x\in \overline{\Omega},
\end{array}\right.
\tag{3.2}
\end{align*}
where $\Omega$ is a bounded, connected open interval in $\mathbb{R}$.
We call a $T$-periodic function $\bar{u}\in C^{1,0}(\mathbb{R}\times \overline{\Omega})$ an upper solution of $(\ref{3.2})$ if $\bar{u}$
satisfies
\begin{align*}
\begin{array}{rl}
\bar{u}_{t}\geq d_{1}[\int_{\Omega}J(x-y)\bar{u}(t,y)dy-\bar{u}(t,x)]+\bar{u}(a(t)-\bar{u}),
\quad (t,x)\in \mathbb{R}\times \overline{\Omega},
\end{array}
\end{align*}
where $C^{1,0}(\mathbb{R}\times \overline{\Omega})$ denotes the class of functions that
are $C^{1}$ in $t$ and continuous in $x$. The lower solution can be defined by reversing the inequality.
For convenience, we define the space $\mathcal{X}_{\Omega}, \mathcal{X}_{\Omega}^{+}, \mathcal{X}_{\Omega}^{++}$
as follows:
\begin{align*}
\begin{array}{rl}
&\mathcal{X}_{\Omega}=\Big\{\phi\in C^{1,0}(\mathbb{R}\times \overline{\Omega}):
\phi(t+T,x)=\phi(t,x), ~(t,x)\in \mathbb{R}\times \overline{\Omega}\Big\},\\[5pt]
&\mathcal{X}_{\Omega}^{+}=\Big\{\phi\in \mathcal{X}_{\Omega}:
\phi(t,x)\geq 0, ~(t,x)\in \mathbb{R}\times \overline{\Omega}\Big\},\\[5pt]
&\mathcal{X}_{\Omega}^{++}=\Big\{\phi\in \mathcal{X}_{\Omega}:
\phi(t,x)> 0, ~(t,x)\in \mathbb{R}\times \overline{\Omega}\Big\}.
\end{array}
\end{align*}

\begin{lemma}\label{l3-3}
Let $\underline{u}\in \mathcal{X}_{\Omega}^{+}, \bar{u}\in \mathcal{X}_{\Omega}^{++}$ be a lower and an upper solution to $(\ref{3.2})$, respectively. Then
$\underline{u}\leq \bar{u}$ in $\mathbb{R}\times \overline{\Omega}$.
\end{lemma}

\begin{proof} The proof follows some ideas of Section 6.3 in \cite{cov10}, where the nonlocal stationary problem was considered.
Define
\begin{align*}
\gamma^{*}:=\inf\{\gamma>0: ~\gamma \bar{u}\geq\underline{u}~\mbox{in}~\mathbb{R}\times \overline{\Omega}\}.
\end{align*}
We claim that $\gamma^{*}\leq 1$. In fact, assume by contradiction that
$\gamma^{*}>1$. Then we have
\begin{align*}\label{3.3}
\begin{array}{rl}
&(\gamma^{*} \bar{u})_{t}-d_{1}[\int_{\Omega}J(x-y)\gamma^{*} \bar{u}(t,y)dy-\gamma^{*} \bar{u}(t,x)]-\gamma^{*} \bar{u}(a(t)-\gamma^{*} \bar{u})\\[5pt]
&\geq \gamma^{*} \bar{u}(a(t)-\bar{u})-\gamma^{*} \bar{u}(a(t)-\gamma^{*} \bar{u})
=\gamma^{*}(\gamma^{*}-1)\bar{u}^{2}>0.
\end{array}
\tag{3.3}
\end{align*}
Since $[0,T]\times \overline{\Omega}$ is compact, we know that $\gamma^{*}$ is attainable, i.e., there exists $(t_{0},x_{0})\in [0,T]\times \overline{\Omega}$ such that $\gamma \bar{u}(t_{0},x_{0})=\underline{u}(t_{0},x_{0})$.

$(i)$ If $(t_{0},x_{0})\in (0,T)\times \overline{\Omega}$, then
$\partial_{t}(\gamma^{*} \bar{u}-\underline{u})(t_{0},x_{0})=0$, since $(t_{0},x_{0})$ is a minimum point
of $\gamma^{*} \bar{u}-\underline{u}$.

$(ii)$ If $(t_{0},x_{0})\in \{0,T\}\times \overline{\Omega}$, by the $T$-periodicity and $C^{1}$-smoothness of $\bar{u}, \underline{u}$ in $t$, we can also deduce $\partial_{t}(\gamma^{*} \bar{u}-\underline{u})(t_{0},x_{0})=0$.

Thus, there holds
\begin{align*}
\begin{array}{rl}
&(\gamma^{*} \bar{u})_{t}(t_{0},x_{0})-d_{1}[\int_{\Omega}J(x_{0}-y)\gamma^{*} \bar{u}(t_{0},y)dy-\gamma^{*} \bar{u}(t_{0},x_{0})]\\[5pt]
&\quad-\gamma^{*} \bar{u}(t_{0},x_{0})(a(t_{0})-\gamma^{*} \bar{u}(t_{0},x_{0}))\\[5pt]
&=\underline{u}_{t}(t_{0},x_{0})-d_{1}[\int_{\Omega}J(x_{0}-y)\gamma^{*} \bar{u}(t_{0},y)dy-\underline{u}(t_{0},x_{0})]\\[5pt]
&\quad-\underline{u}(t_{0},x_{0})(a(t_{0})-\underline{u}(t_{0},x_{0}))\\[5pt]
&\leq d_{1}\int_{\Omega}J(x_{0}-y)[\underline{u}(t_{0},y)-\gamma^{*} \bar{u}(t_{0},y)]dy
\leq 0,
\end{array}
\end{align*}
which contradicts with $(\ref{3.3})$. Therefore, the claim is true. It follows that
$\underline{u}\leq \bar{u}$ in $[0,T]\times \overline{\Omega}$.
\end{proof}

\subsection{Some eigenvalue problems}

In this subsection, we mainly study some eigenvalue problems and analyze the properties
of their principle eigenvalue. Hereafter, we always assume $\Omega$ be a bounded, connected open interval in $\mathbb{R}$
and $|\Omega|$ be its length.

For $(t,x)\in \mathbb{R}\times \overline{\Omega}$, we consider the following operator
\begin{align*}\label{3.4}
\begin{array}{rl}
-(L_{\Omega}+a)[\phi](t,x)=
\phi_{t}(t,x)-d_{1}[\int_{\Omega}J(x-y)\phi(t,y)dy-\phi(t,x)]-a(t)\phi(t,x),
\end{array}
\tag{3.4}
\end{align*}
where
$a\in C_{T}(\mathbb{R}):=\{a\in C(\mathbb{R}): a(t+T)=a(t)>0, \forall t\in \mathbb{R}\}$.

We define
\begin{align*}
\lambda_{1}(-(L_{\Omega}+a))=\inf\Big\{\mathfrak{R}\lambda:~\lambda\in\sigma(-(L_{\Omega}+a))\Big\},
\end{align*}
where $\sigma(-(L_{\Omega}+a))$ is the spectrum of $-(L_{\Omega}+a)$.
By Theorem \textbf{A}$(1)$ in \cite{svo19}, we know that $\lambda_{1}(-(L_{\Omega}+a))$
is the principle eigenvalue of $-(L_{\Omega}+a)$, which means that there exists an eigenfunction
$\phi\in \mathcal{X}_{\Omega}^{++}$ such that
\begin{align*}
-(L_{\Omega}+a)[\phi](t,x)=\lambda_{1}(-(L_{\Omega}+a))\phi.
\end{align*}

\begin{lemma}\label{l3.4} (see Theorem \textbf{B} in \cite{svo19})
Assume that $J$ satisfies $(\ref{J})$ and $a\in C_{T}(\mathbb{R})$.
Let $u(t,x;u_{0})$ be a solution of
\begin{align*}
\left\{\begin{array}{l}
u_t=d_{1}[\int_{\Omega}J(x-y)u(t,y)dy-u(t,x)]+u(a(t)-u),
\quad t>0, x\in \overline{\Omega},\\[3pt]
u(0,x)=u_{0}(x), \quad x\in \overline{\Omega},
\end{array}\right.
\end{align*}
where $u_{0}\in C(\overline{\Omega})$ is non-negative and not identically zero.
The following statements hold:

$(i)$ If $\lambda_{1}(-(L_{\Omega}+a))<0$, then the equation
\begin{align*}\label{3.5}
\begin{array}{rl}
u_t=d_{1}[\int_{\Omega}J(x-y)u(t,y)dy-u(t,x)]+u(a(t)-u),
\quad t\in \mathbb{R}, x\in \overline{\Omega}
\end{array}
\tag{3.5}
\end{align*}
admits a unique solution $u^{*}\in \mathcal{X}_{\Omega}^{++}$, and there holds
\begin{align*}
\|u(t,\cdot;u_{0})-u^{*}(t,\cdot)\|_{C(\overline{\Omega})}\rightarrow 0
\quad \mbox{as}~t\rightarrow \infty,
\end{align*}

$(ii)$ If $\lambda_{1}(-(L_{\Omega}+a))>0$, then the equation
$(\ref{3.5})$ admits no solution in $\mathcal{X}_{\Omega}^{+}\setminus \{0\}$
and there holds
\begin{align*}
\|u(t,\cdot;u_{0})\|_{C(\overline{\Omega})}\rightarrow 0
\quad \mbox{as}~t\rightarrow \infty.
\end{align*}
\end{lemma}

\begin{remark}\label{r3.5} For the case $\lambda_{1}(-(L_{\Omega}+a))=0$,
$(\ref{3.5})$ has been shown in \cite{svo19} to admit no solution in $\mathcal{X}_{\Omega}^{+}\setminus \{0\}$, but the global dynamics is not provided. Since $a(t)$ is independent of spatial variable, we can also get $\|u(t,\cdot;u_{0})\|_{C(\overline{\Omega})}\rightarrow 0$, more details can be seen in the proof of Theorem $\ref{t4.4}$.
\end{remark}

In what follows, we present some further properties of $\lambda_1$.

\begin{lemma}\label{l3.6}
Let $J$ satisfies $(\ref{J})$ and $a\in C_{T}(\mathbb{R})$. Then

$(i)$ $\lambda_{1}(-(L_{\Omega}+a))$ is strictly decreasing and continuous in
$|\Omega|$;

$(ii)$ $\lim_{|\Omega|\rightarrow +\infty}\lambda_{1}(-(L_{\Omega}+a))=-a_{T}$,
where $a_{T}=\frac{1}{T}\int_{0}^{T}a(t)dt$;

$(iii)$ $\lim_{|\Omega|\rightarrow 0}\lambda_{1}(-(L_{\Omega}+a))=d_{1}-a_{T}$.
\end{lemma}

\begin{proof} Let $\phi\in \mathcal{X}_{\Omega}^{++}$ be an eigenfunction of
$-(L_{\Omega}+a)$ associated with the principle eigenvalue $\lambda_{1}(-(L_{\Omega}+a))$.
We define
\begin{align*}
\psi(t,x)=e^{-\int_{0}^{t}(a(s)-a_{T})ds}\phi(t,x), \quad\forall (t,x)\in \mathbb{R}\times \overline{\Omega}.
\end{align*}
It is easy to check that $\psi\in \mathcal{X}_{\Omega}^{++}$.

Multiplying the equation $-(L_{\Omega}+a)[\phi]=\lambda_{1}(-(L_{\Omega}+a))\phi$
by the function $t\mapsto e^{-\int_{0}^{t}(a(s)-a_{T})ds}$, we have
\begin{align*}
\begin{array}{rl}
-\psi_{t}(t,x)+d_{1}[\int_{\Omega}J(x-y)\psi(t,y)dy-\psi(t,x)]
+a_{T}\psi(t,x)+\lambda_{1}(-(L_{\Omega}+a))\psi(t,x)=0
\end{array}
\end{align*}
for $(t,x)\in \mathbb{R}\times \overline{\Omega}$. Taking $\psi_{T}(x)=\frac{1}{T}\int_{0}^{T}\psi(t,x)dt$ for $x\in \overline{\Omega}$,
and integrating the above equation over $[0,T]$ with respect to $t$, we have
\begin{align*}
\begin{array}{rl}
d_{1}[\int_{\Omega}J(x-y)\psi_{T}(y)dy-\psi_{T}(x)]
+a_{T}\psi_{T}(x)+\lambda_{1}(-(L_{\Omega}+a))\psi_{T}(x)=0,
\quad x\in \overline{\Omega}.
\end{array}
\end{align*}
That is, $\lambda_{1}(-(L_{\Omega}+a))$ is the principle eigenvalue of the following nonlocal operator
$-(\mathcal{L}_{\Omega}+a_{T}): C(\overline{\Omega})\rightarrow C(\overline{\Omega})$ defined by
\begin{align*}\label{3.6}
\begin{array}{rl}
-(\mathcal{L}_{\Omega}+a_{T})[\omega](x)
:=-d_{1}[\int_{\Omega}J(x-y)\omega(y)dy-\omega(x)]
-a_{T}\omega(x)
\end{array}
\tag{3.6}
\end{align*}
with an eigenfunction $\psi_{T}\in \mathcal{X}_{\Omega}^{++}$. Denote by $\lambda_{1}(-(\mathcal{L}_{\Omega}+a_{T}))$ the principle eigenvalue of $-(\mathcal{L}_{\Omega}+a_{T})$,
then we have
\begin{align*}\label{3.7}
\lambda_{1}(-(L_{\Omega}+a))=\lambda_{1}(-(\mathcal{L}_{\Omega}+a_{T})).
\tag{3.7}
\end{align*}

Without loss of generality, we assume that $\Omega=(l_{1},l_{2})$.
According to Proposition 3.4 in \cite{cdll19}, we know the following results hold:

$(i)$ $\lambda_{1}(-(\mathcal{L}_{\Omega}+a_{T}))$ is strictly increasing and continuous in
$|\Omega|=l_{2}-l_{1}$;

$(ii)$ $\lim_{l_{2}-l_{1}\rightarrow +\infty}\lambda_{1}(-(\mathcal{L}_{\Omega}+a_{T}))=a_{T}$;

$(iii)$ $\lim_{l_{2}-l_{1}\rightarrow 0}\lambda_{1}(-(\mathcal{L}_{\Omega}+a_{T}))=a_{T}-d_{1}$.\\
Combining the above conclusions and $(\ref{3.7})$, we can get the desired results.
\end{proof}

Now, we consider another periodic-parabolic eigenvalue problem
\begin{align*}\label{3.8}
\left\{\begin{array}{l}
-(\tilde{L}_{\Omega}+c)[\varphi](t,x)=\varphi_{t}-d_{2}[\tau\varphi_{xx}+(1-\tau)(\int_{\Omega}J(x-y)\varphi(t,y)dy-\varphi)]-c(t)\varphi\\[3pt]
\qquad\qquad\qquad\qquad=\lambda \varphi \quad \mbox{in}~[0,T]\times\Omega,\\[3pt]
\varphi(t,x)=0 \quad \mbox{on}~[0,T]\times\partial\Omega, \\[3pt]
\varphi(0,x)=\varphi(T,x) \quad \mbox{in}~\Omega.
\end{array}\right.
\tag{3.8}
\end{align*}
Define a linear nonlocal operator $\mathcal{K}$ as follows
\begin{align*}
\begin{array}{rl}
(\mathcal{K}\varphi)(t,x):=
\int_{\Omega}J(x-y)\varphi(t,y)dy-\varphi(t,x).
\end{array}
\end{align*}
For any given $0<\tau\leq 1$, we can check that $\{\mathcal{A}(t): 0\leq t\leq T\}:=\{-d_{2}[\tau\partial_{x}^{2}+(1-\tau)\mathcal{K}]-c(t)I: 0\leq t\leq T\}$ satisfy the hypotheses
(11.5) in \cite{hess91}. As showed in Section II.14 of \cite{hess91},
based on the Krein-Rutman theorem, we can prove that $(\ref{3.8})$ admits a principle eigenvalue
$\lambda_{1}(-(\tilde{L}_{\Omega}+c))$ with principle eigenfunction $\varphi$.\\

For later applications, we give the following lemma.

\begin{lemma}\label{l3.7} Let $J$ satisfies $(\ref{J})$ and $c\in C_{T}(\mathbb{R})$. Then

$(i)$ $\lambda_{1}(-(\tilde{L}_{\Omega}+c))$ is a strictly decreasing continuous function in
$|\Omega|$. Moreover, $\lim_{|\Omega|\rightarrow 0}\lambda_{1}(-(\tilde{L}_{\Omega}+c))=+\infty$,
$\lim_{|\Omega|\rightarrow +\infty}\lambda_{1}(-(\tilde{L}_{\Omega}+c))=-c_{T}$. Then  $\lambda_{1}(-(\tilde{L}_{\Omega}+c))=0$ has a unique root $|\Omega|=h^{*}$;

$(ii)$ if $\lambda_{1}(-(\tilde{L}_{\Omega}+c))<0$, then the problem
\begin{align*}
\left\{\begin{array}{l}
\varphi_{t}-d_{2}[\tau\varphi_{xx}+(1-\tau)(\int_{\Omega}J(x-y)\varphi(t,y)dy-\varphi(t,x))]\\[3pt]
\quad=\varphi(c(t)-\varphi) \quad \mbox{in}~(0,\infty)\times\Omega,\\[3pt]
\varphi(t,x)=0 \quad \mbox{on}~(0,\infty)\times\partial\Omega.
\end{array}\right.
\end{align*}
admits a unique positive $T$-periodic solution $\varphi^{*}$,
and $\varphi^{*}$ is globally asymptotically stable.
\end{lemma}

\begin{proof}
$(i)$ Let $\varphi$ be an eigenfunction of
$(\ref{3.8})$ associated with the principle eigenvalue $\lambda_{1}(-(\tilde{L}_{\Omega}+c))$. Define
\begin{align*}
\psi(t,x)=e^{-\int_{0}^{t}(c(s)-c_{T})ds}\varphi(t,x), \quad\forall (t,x)\in \mathbb{R}\times \overline{\Omega}.
\end{align*}
Similar as the proof of Lemma $\ref{l3.6}$, $\lambda_{1}(-(\tilde{L}_{\Omega}+c))$ is the principal eigenvalue of the following elliptic-type problem
\begin{align*}\label{3.9}
\left\{\begin{array}{l}
-(\tilde{\mathcal{L}}_{\Omega}+c_{T})[\omega]=-d_{2}[\tau\omega_{xx}+(1-\tau)(\int_{\Omega}J(x-y)\omega(y)dy-\omega(x))]-c_{T}\omega\\[3pt]
\qquad\qquad\qquad\quad=\lambda \omega \quad \mbox{in}~\Omega,\\[3pt]
\omega(x)=0 \quad \mbox{on}~\partial\Omega
\end{array}\right.
\tag{3.9}
\end{align*}
with an eigenfunction $\omega(x)=\frac{1}{T}\int_{0}^{T}\psi(t,x)dt$.
Denote by $\lambda_{1}(-(\tilde{\mathcal{L}}_{\Omega}+c_{T}))$ the principle eigenvalue of $(\ref{3.9})$,
then we have
\begin{align*}\label{3.10}
\lambda_{1}(-(\tilde{L}_{\Omega}+c))=\lambda_{1}(-(\tilde{\mathcal{L}}_{\Omega}+c_{T})).
\tag{3.10}
\end{align*}
The continuity of $\lambda_{1}(-(\tilde{\mathcal{L}}_{\Omega}+c_{T}))$ with respect to $|\Omega|$ can be obtained by using a simple re-scaling argument of the spatial variable $x$.
Note that $\lambda_{1}(-(\tilde{\mathcal{L}}_{\Omega}+c_{T}))$ can be expressed in a variational formulation
\begin{align*}
\begin{array}{rl}
&\lambda_{1}(-(\tilde{\mathcal{L}}_{\Omega}+c_{T}))\\[3pt]
&=\inf_{0\not\equiv \omega\in H_{0}^{1}(\Omega)}
\frac{d_{2}\tau\int_{\Omega}\omega_{x}^{2}(x)dx-d_{2}(1-\tau)\int_{\Omega}\int_{\Omega}J(x-y)\omega(y)\omega(x)dydx}
{\int_{\Omega}\omega^{2}(x)dx}+[d_{2}(1-\tau)-c_{T}].
\end{array}
\end{align*}
By the zero extension of principle eigenfunction,
we can get the monotonicity of $\lambda_{1}(\tilde{\mathcal{L}}_{\Omega}+c_{T})$
from the variational formulation of principle eigenvalue.

Next, we consider the limits of $\lambda_{1}(-(\tilde{L}_{\Omega}+c))$ as $|\Omega|\rightarrow 0$
and $|\Omega|\rightarrow\infty$.
Without loss of generality, we may assume that $\Omega=(0,l)$.
Since
\begin{align*}
\begin{array}{rl}
\int_{0}^{l}\int_{0}^{l}J(x-y)\omega(y)\omega(x)dydx
\leq \int_{0}^{l}\int_{0}^{l}J(x-y)\frac{\omega^{2}(y)+\omega^{2}(x)}{2}dydx
\leq \int_{0}^{l}\omega^{2}(x)dx,
\end{array}
\end{align*}
we have
\begin{align*}
\begin{array}{rl}
\lambda_{1}(-(\tilde{\mathcal{L}}_{(0,l)}+c_{T}))
\geq \inf_{0\not\equiv \omega\in H_{0}^{1}((0,l))}
\frac{d_{2}\tau\int_{0}^{l}\omega_{x}^{2}(x)dx}
{\int_{0}^{l}\omega^{2}(x)dx}-c_{T}.
\end{array}
\end{align*}
By the fact that
\begin{align*}
\begin{array}{rl}
\inf_{0\not\equiv \omega\in H_{0}^{1}((0,l))}
\frac{\int_{0}^{l}\omega_{x}^{2}(x)dx}
{\int_{0}^{l}\omega^{2}(x)dx}
=\frac{\pi^{2}}{4l^{2}},
\end{array}
\end{align*}
we know
\begin{align*}\label{3.11}
\begin{array}{rl}
\lim_{l\rightarrow 0}\lambda_{1}(-(\tilde{\mathcal{L}}_{(0,l)}+c_{T}))=+\infty
\end{array}
\tag{3.11}
\end{align*}
and
\begin{align*}\label{3.12}
\begin{array}{rl}
\liminf_{l\rightarrow +\infty}\lambda_{1}(-(\tilde{\mathcal{L}}_{(0,l)}+c_{T}))\geq -c_{T}.
\end{array}
\tag{3.12}
\end{align*}

On the other hand, by $(\ref{J})$, for any fixed $0<\varepsilon\ll 1$, there exists $L=L(\varepsilon)>0$
such that
\begin{align*}
\begin{array}{rl}
\int_{-L}^{L}J(x)dx>1-\varepsilon.
\end{array}
\end{align*}
For any large $l>3L$, we choose the test function $\varphi_{\varepsilon}(x)$ defined as follows
\begin{align*}
\varphi_{\varepsilon}(x)
=\left\{\begin{array}{l}
\frac{x}{\varepsilon}, \quad x\in[0,\varepsilon],\\[3pt]
1, \quad x\in[\varepsilon, l-\varepsilon],\\[3pt]
\frac{l-x}{\varepsilon}, \quad x\in[l-\varepsilon,l].
\end{array}\right.
\end{align*}
It is easy to check that $\varphi_{\varepsilon}\in H_{0}^{1}((0,l))$ and satisfies
$\int_{0}^{l}\varphi_{\varepsilon}^{2}(x)dx=l-\frac{4}{3}\varepsilon$
and $\int_{0}^{l}(\partial_{x}\varphi_{\varepsilon})^{2}(x)dx=\frac{2}{\varepsilon}$.
Thus,
\begin{align*}
\begin{array}{rl}
&\lambda_{1}(-(\tilde{\mathcal{L}}_{(0,l)}+c_{T}))\\[3pt]
&\leq \frac{d_{2}\tau\int_{0}^{l}(\partial_{x}\varphi_{\varepsilon})^{2}(x)dx
-d_{2}(1-\tau)\int_{0}^{l}\int_{0}^{l}J(x-y)\varphi_{\varepsilon}(y)\varphi_{\varepsilon}(x)dydx}
{\int_{0}^{l}\varphi_{\varepsilon}^{2}(x)dx}+[d_{2}(1-\tau)-c_{T}]\\[3pt]
&\leq \frac{\frac{2d_{2}\tau}{\varepsilon}
-d_{2}(1-\tau)\int_{L+\varepsilon}^{l-L-\varepsilon}\int_{\varepsilon}^{l-\varepsilon}J(x-y)dydx}
{l-\frac{4}{3}\varepsilon}+[d_{2}(1-\tau)-c_{T}]\\[3pt]
&\leq \frac{\frac{2d_{2}\tau}{\varepsilon}
-d_{2}(1-\tau)\int_{L+\varepsilon}^{l-L-\varepsilon}\int_{-L}^{L}J(\xi)d\xi dx}
{l-\frac{4}{3}\varepsilon}+[d_{2}(1-\tau)-c_{T}]\\[3pt]
&\leq \frac{\frac{2d_{2}\tau}{\varepsilon}
-d_{2}(1-\tau)(l-2L-2\varepsilon)(1-\varepsilon)}
{l-\frac{4}{3}\varepsilon}+[d_{2}(1-\tau)-c_{T}]\\[3pt]
&\rightarrow -d_{2}(1-\tau)(1-\varepsilon)+[d_{2}(1-\tau)-c_{T}] \quad
\mbox{as}~l\rightarrow +\infty.
\end{array}
\end{align*}
Since $\varepsilon$ is arbitrary, it follows that
\begin{align*}
\begin{array}{rl}
\limsup_{l\rightarrow +\infty}\lambda_{1}(-(\tilde{\mathcal{L}}_{(0,l)}+c_{T}))\leq -c_{T},
\end{array}
\end{align*}
which together with $(\ref{3.12})$ imply that
\begin{align*}\label{3.13}
\begin{array}{rl}
\lim_{l\rightarrow +\infty}\lambda_{1}(-(\tilde{\mathcal{L}}_{(0,l)}+c_{T}))=-c_{T}.
\end{array}
\tag{3.13}
\end{align*}
From $(\ref{3.10})$, $(\ref{3.11})$ and $(\ref{3.13})$, we know that $\lambda_{1}(-(\tilde{L}_{(0,l)}+c))=0$ has a unique root.

$(ii)$ the proof is similar as that of Lemma 3.3 in \cite{w143} or Theorem 28.1 in \cite{hess91}, we omit
the details.
\end{proof}

\section{Spreading and vanishing for problem $(\ref{1.1})$}

In this section, we investigate the dynamics of problem $(\ref{1.1})$, including the spreading-vanishing dichotomy and some sufficient conditions for spreading and vanishing.
In view of $(\ref{2.2})$, we see that the free boundaries $h(t),-g(t)$ are strictly increasing functions with respect to time $t$. Thus, $h_{\infty}:=\lim_{t\rightarrow \infty}h(t)$ and
$g_{\infty}:=\lim_{t\rightarrow \infty}g(t)$ are well-defined. Clearly, $h_{\infty},-g_{\infty}\leq +\infty$.

By similar argument as the proof of Proposition 3.1 in \cite{w14} with minor modifications, we have the following result.

\begin{lemma}\label{l4.1}
Let $d$, $\mu$ and $h^{0}$ be positive constants and $C\in \mathbb{R}$. Assume that $\varphi_{0}\in C^{2}([-h^{0}, h^{0}])$ satisfies $\varphi_{0}(-h^{0})=\varphi_{0}(h^{0})=0$ and $\varphi_{0}>0$ in $(-h^{0}, h^{0})$. Let
$(g, h)\in [C^{1+\frac{\alpha}{2}}[0, \infty)]^{2}$, $\varphi\in C^{1+\frac{\alpha}{2}, 2+\alpha}((0, \infty)\times(g(t), h(t)))$ for some $\alpha\in (0,1)$ and
satisfy $g(t)<0$, $h(t)>0$, $\varphi(t,x)>0$ for all $t\geq 0$ and $g(t)<x<h(t)$. We further suppose that
$\lim_{t\rightarrow \infty}g(t)>-\infty$, $\lim_{t\rightarrow \infty}h(t)<\infty$,
$\lim_{t\rightarrow \infty}g^{\prime}(t)=\lim_{t\rightarrow \infty}h^{\prime}(t)=0$
and there exists a constant $K>0$ such that $\|\varphi\|_{C^{1}[g(t), h(t)]}\leq K$ for $t>1$.
If $(\varphi, g, h)$ satisfies
\begin{align*}
\left\{\begin{array}{l}
\varphi_t-d\varphi_{xx}\geq C\varphi,\quad t>0,~ g(t)<x<h(t),\\[3pt]
\varphi=0,\quad t\geq0, ~x=g(t)~\mbox{or}~ x=h(t),\\[3pt]
g^{\prime}(t)\leq-\mu\varphi_{x}(t, g(t)),~h^{\prime}(t)\geq-\mu\varphi_{x}(t, h(t)),
\quad t>0,\\[3pt]
g(0)=-h^{0},~ h(0)=h^{0}, \\[3pt]
\varphi(0,x)=\varphi_{0}(x), \quad -h^{0}<x<h^{0},
\end{array}\right.
\end{align*}
then $\lim_{t\rightarrow \infty}\max_{g(t)\leq x\leq h(t)}\varphi(t,x)=0$.
\end{lemma}

The next lemma provides an estimate for $v$. The proof is a simple modification of that for Lemma 3.2 in \cite{waw181}, so we omit it here.

\begin{lemma}\label{l4.2} Let $(u,v,g,h)$ be the unique global solution of $(\ref{1.1})$ and
$h_{\infty}-g_{\infty}<\infty$. Then
there exists $C>0$ such that
\begin{align*}\label{4.1}
\begin{array}{rl}
\|v\|_{C^{\frac{1+\alpha}{2},1+\alpha}(D_{\infty})}\leq C,
\quad \mbox{where}~D_{\infty}:=[0,\infty)\times[g(t),h(t)]
\end{array}
\tag{4.1}
\end{align*}
and hence
\begin{align*}\label{4.2}
\begin{array}{rl}
\|v_{x}(t,g(t))\|_{C^{\frac{\alpha}{2}}(\overline{\mathbb{R}}_{+})}
+\|v_{x}(t,h(t))\|_{C^{\frac{\alpha}{2}}(\overline{\mathbb{R}}_{+})}
\leq C.
\end{array}
\tag{4.2}
\end{align*}
\end{lemma}

\begin{lemma}\label{l4.3} If $h_{\infty}-g_{\infty}<\infty$,
then $\lim_{t\rightarrow \infty}g^{\prime}(t)=\lim_{t\rightarrow \infty}h^{\prime}(t)=0$.
\end{lemma}

\begin{proof} It is easy to see that $-\infty<g_{\infty}<h_{\infty}<\infty$.
From $(\ref{2.2})$, we can deduce that $g^{\prime}(t)$ and $h^{\prime}(t)$ defined in $(\ref{1.1})$ are bounded. Let
\begin{align*}
\begin{array}{l}
\varphi_{1}(t)=v_{x}(t,h(t)),
~\varphi_{2}(t)=\int_{g(t)}^{h(t)}\int_{h(t)}^{\infty}J(x-y)u(t,x)dydx,\\[3pt]
\varphi_{3}(t)=\int_{g(t)}^{h(t)}\int_{h(t)}^{\infty}J(x-y)v(t,x)dydx.
\end{array}
\end{align*}
By $(\ref{4.2})$, we get
$|\varphi_{1}(t)-\varphi_{1}(s)|\leq C_{1}|t-s|^{\frac{\alpha}{2}}$ for any $t,s>0$.
For $\varphi_{2}$, assume $t>s$, we have $h(t)>h(s)$,
$g(t)<g(s)$ and then
\begin{align*}
\begin{array}{rl}
&\varphi_{2}(t)-\varphi_{2}(s)\\[3pt]
&=\int_{g(t)}^{h(t)}\int_{h(t)}^{\infty}J(x-y)u(t,x)dydx
-\int_{g(s)}^{h(s)}\int_{h(s)}^{\infty}J(x-y)u(s,x)dydx\\[3pt]
&=\int_{g(s)}^{h(s)}\int_{h(t)}^{\infty}J(x-y)[u(t,x)-u(s,x)]dydx
+\int_{g(t)}^{g(s)}\int_{h(t)}^{\infty}J(x-y)u(t,x)dydx\\[3pt]
&\quad+\int_{h(s)}^{h(t)}\int_{h(t)}^{\infty}J(x-y)u(t,x)dydx
-\int_{g(s)}^{h(s)}\int_{h(s)}^{h(t)}J(x-y)u(s,x)dydx\\[3pt]
&\leq \|\partial_{t}u\|_{L^{\infty}(D_{\infty})}(t-s)(h(s)-g(s))
+\|u\|_{L^{\infty}(D_{\infty})}(g(s)-g(t))\\[3pt]
&\quad+2\|u\|_{L^{\infty}(D_{\infty})}(h(t)-h(s))\\[3pt]
&\leq C_{2}(t-s),
\end{array}
\end{align*}
where $\|\partial_{t}u\|_{L^{\infty}(D_{\infty})}$ is obtained by
the first equation in $(\ref{1.1})$ and the bound of $u$.
Thus,
\begin{align*}
|\varphi_{2}(t)-\varphi_{2}(s)|
\leq C_{2}|t-s|.
\end{align*}

For $\varphi_{3}$, it follows from $(\ref{4.1})$ that $|v(t,x)-v(s,x)|\leq C|t-s|^{\frac{1+\alpha}{2}}$
for $x\in [g(t),h(t)]$. Similar to $\varphi_{2}$, we can prove that
\begin{align*}
|\varphi_{3}(t)-\varphi_{3}(s)|
\leq C_{3}|t-s|.
\end{align*}
Therefore, $h^{\prime}(t)=-\mu \varphi_{1}+\rho_{1}\varphi_{2}+\rho_{2}\varphi_{3}$
is uniformly continuous in $[0,\infty)$.
From $\lim_{t\rightarrow \infty}h(t)=h_{\infty}<\infty$, we know
$\lim_{t\rightarrow \infty}h^{\prime}(t)=0$. Similarly, we can show
$\lim_{t\rightarrow \infty}g^{\prime}(t)=0$.
\end{proof}

\begin{theorem}\label{t4.4} If
$h_{\infty}-g_{\infty}<\infty$, then the solution $(u,v,g,h)$ of $(\ref{1.1})$ satisfies
\begin{align*}
\lim_{t\rightarrow \infty}\|u(t,\cdot)\|_{C([g(t),h(t)])}=
\lim_{t\rightarrow \infty}\|v(t,\cdot)\|_{C([g(t),h(t)])}=0.
\end{align*}
\end{theorem}

\begin{proof} Since $J\geq 0$ and $v>0$, from the second equation in $(\ref{1.1})$,
there exists a constant $C>0$ such that
\begin{align*}
\partial_{t}v-d_{2}\tau\partial_{x}^{2}v\geq Cv.
\end{align*}
According to Lemma $\ref{l4.1}$, we get
\begin{align*}
\lim_{t\rightarrow \infty}\|v(t,\cdot)\|_{C([g(t),h(t)])}=0.
\end{align*}
We claim that
\begin{align*}\label{4.3}
\lambda_{1}(-(L_{(g_{\infty},h_{\infty})}+a))\geq 0,
\tag{4.3}
\end{align*}
where $-(L_{(g_{\infty},h_{\infty})}+a)$ is defined in $(\ref{3.4})$.
Assume on the contrary that $\lambda_{1}(-(L_{(g_{\infty},h_{\infty})}+a))<0$.
For convenient, for any $\varepsilon>0$ we define
$h_{\infty}^{\pm \varepsilon}:=h_{\infty}\pm \varepsilon$,
$g_{\infty}^{\pm \varepsilon}:=g_{\infty}\pm \varepsilon$.
Thus, there exists $\varepsilon_{1}>0$ such that $\lambda_{1}(-(L_{(g_{\infty}^{+\varepsilon},h_{\infty}^{-\varepsilon})}+a(t)-b(t)\varepsilon))>0$
for all $\varepsilon\in (0, \varepsilon_{1})$. For such $\varepsilon>0$, we can find $T_{\varepsilon}>0$ such that, for $t>T_{\varepsilon}$,
\begin{align*}
h(t)>h_{\infty}^{-\varepsilon},~g(t)<g_{\infty}^{+\varepsilon},~
\|v(t,\cdot)\|_{C([g(t),h(t)])}<\varepsilon.
\end{align*}
Then $u$ satisfies
\begin{align*}
\left\{\begin{array}{l}
u_t\geq d_{1}\int_{g_{\infty}^{+\varepsilon}}^{h_{\infty}^{-\varepsilon}}
J(x-y)u(t,y)dy-d_{1}u+u(a(t)-u-b(t)\varepsilon),\quad t>T_{\varepsilon},~
x\in [g_{\infty}^{+\varepsilon}, h_{\infty}^{-\varepsilon}],\\[5pt]
u(T_{\varepsilon},x)=u(T_{\varepsilon},x), \quad x\in [g_{\infty}^{+\varepsilon}, h_{\infty}^{-\varepsilon}].
\end{array}\right.
\end{align*}
Consider the following problem
\begin{align*}\label{4.4}
\left\{\begin{array}{l}
\phi_t= d_{1}\int_{g_{\infty}^{+\varepsilon}}^{h_{\infty}^{-\varepsilon}}
J(x-y)\phi(t,y)dy-d_{1}\phi+\phi(a(t)-\phi-b(t)\varepsilon),\quad t>T_{\varepsilon},~
x\in [g_{\infty}^{+\varepsilon}, h_{\infty}^{-\varepsilon}],\\[5pt]
\phi(T_{\varepsilon},x)=u(T_{\varepsilon},x), \quad x\in [g_{\infty}^{+\varepsilon}, h_{\infty}^{-\varepsilon}].
\end{array}\right.
\tag{4.4}
\end{align*}
Since $\lambda_{1}(-(L_{(g_{\infty}^{+\varepsilon},h_{\infty}^{-\varepsilon})}+a(t)-b(t)\varepsilon))<0$,
by Lemma $\ref{l3.4}$$(i)$ we know that the solution $\phi_{\varepsilon}(t,x)$ of problem $(\ref{4.4})$ converges to $\phi^{*}_{\varepsilon}(t,x)$ uniformly in $[g_{\infty}^{+\varepsilon}, h_{\infty}^{-\varepsilon}]$ as $t\rightarrow \infty$, where $\phi^{*}_{\varepsilon}(t,x)\in \mathcal{X}_{\varepsilon}^{++}$ is the unique periodic solution of
\begin{align*}
\phi_t=d_{1}\int_{g_{\infty}^{+\varepsilon}}^{h_{\infty}^{-\varepsilon}}
J(x-y)\phi(t,y)dy-d_{1}\phi+\phi(a(t)-\phi-b(t)\varepsilon),\quad t\in \mathbb{R},~
x\in [g_{\infty}^{+\varepsilon}, h_{\infty}^{-\varepsilon}].
\end{align*}
By Lemma 3.3 in \cite{cdll19} and a simple comparison argument, we get
\begin{align*}
u(t,x)\geq \phi_{\varepsilon}(t,x), \quad
\forall~t>T_{\varepsilon},~x\in [g_{\infty}^{+\varepsilon}, h_{\infty}^{-\varepsilon}].
\end{align*}
Hence, there exist two constants $\tilde{T}_{\varepsilon}>T_{\varepsilon}$ and $C>0$
such that
\begin{align*}
u(t,x)\geq \frac{1}{2}\phi^{*}_{\varepsilon}(t,x)\geq C>0, \quad
\forall~t>\tilde{T}_{\varepsilon},~x\in [g_{\infty}^{+\varepsilon}, h_{\infty}^{-\varepsilon}].
\end{align*}
It follows that, for $0<\varepsilon<\min\{\varepsilon_{1}, \frac{\bar{\epsilon}}{2}\}$ and $t>\tilde{T}_{\varepsilon}$,
\begin{align*}
h^{\prime}(t)
&\geq \rho_{1}\int_{g(t)}^{h(t)}\int_{h(t)}^{\infty}J(x-y)u(t,x)dydx
\geq \rho_{1}\int_{g_{\infty}^{+\varepsilon}}^{h_{\infty}^{-\varepsilon}}\int_{h_{\infty}}^{\infty}J(x-y)u(t,x)dydx\\[5pt] &\geq \rho_{1}\int_{h_{\infty}^{-\frac{\bar{\varepsilon}}{2}}}^{h_{\infty}^{-\varepsilon}}\int_{h_{\infty}}^{h_{\infty}^{-\frac{\bar{\varepsilon}}{2}}}
\eta_{0}\frac{1}{2}\phi^{*}_{\varepsilon}(t,x)dydx
\geq \rho_{1}\int_{h_{\infty}^{-\frac{\bar{\varepsilon}}{2}}}^{h_{\infty}^{-\varepsilon}}\int_{h_{\infty}}^{h_{\infty}^{-\frac{\bar{\varepsilon}}{2}}}
\eta_{0}Cdydx>0,
\end{align*}
which implies that $h_{\infty}=\infty$. It is a contradiction and then $(\ref{4.3})$ holds.

Let $\bar{u}$ be the unique solution of
\begin{align*}
\left\{\begin{array}{l}
\bar{u}_t=d_{1}\int_{g_{\infty}}^{h_{\infty}}
J(x-y)\bar{u}(t,y)dy-d_{1}\bar{u}+\bar{u}(a(t)-\bar{u}),\quad t>0,
x\in [g_{\infty}, h_{\infty}],\\[5pt]
\bar{u}(0,x)=u_{0}(x),~ x\in [-h_{0},h_{0}];
\quad \bar{u}(0,x)=0,
~x\in [g_{\infty}, h_{\infty}]\setminus[-h_{0}, h_{0}].
\end{array}\right.
\end{align*}
Now we prove that $\lim_{t\rightarrow \infty}\bar{u}(t,x)=0$ uniformly
in $[g_{\infty},h_{\infty}]$. Since $(\ref{4.3})$ holds, we divide the discussion
into two cases:

$(i)$ For the case $\lambda_{1}(-(L_{(g_{\infty},h_{\infty})}+a))>0$,
applying Lemma $\ref{l3.4}$$(ii)$ we can get the desired result.

$(ii)$ For the case $\lambda_{1}(-(L_{(g_{\infty},h_{\infty})}+a))=0$,
we define
\begin{align*}
w(t,x)=e^{-\int_{0}^{t}[a(s)-a_{T}]ds}\bar{u}(t,x),
\end{align*}
then $w(t,x)$ satisfies
\begin{align*}
\left\{\begin{array}{l}
w_t=d_{1}\int_{g_{\infty}}^{h_{\infty}}
J(x-y)w(t,y)dy-d_{1}w+w(a_{T}-e^{\int_{0}^{t}[a(s)-a_{T}]ds}w),~t>0,
x\in [g_{\infty}, h_{\infty}],\\[5pt]
w(0,x)=u_{0}(x),~x\in [-h_{0}, h_{0}];
\quad w(0,x)=0,~x\in [g_{\infty}, h_{\infty}]\setminus[-h_{0}, h_{0}].
\end{array}\right.
\end{align*}
For any $t>0$, we can write $t=nT+\tau$ with $\tau\in [0,T)$, and then
\begin{align*}
e^{\int_{0}^{t}[a(s)-a_{T}]ds}
=e^{\int_{0}^{t}[a(s)-a_{T}]ds}
=e^{\int_{0}^{\tau}[a(s)-a_{T}]ds},
\end{align*}
which together with the continuity of $a(t)$ imply that $M_{1}\leq e^{\int_{0}^{t}[a(s)-a_{T}]ds}\leq M_{2}$ for some
positive constants $M_{1}$ and $M_{2}$.
By the comparison principle (Lemma 3.3 in \cite{cdll19}), we know $w(t,x)\leq \bar{w}(t,x)$ with $\bar{w}(t,x)$
be the unique solution of
\begin{align*}
\left\{\begin{array}{l}
\bar{w}_t=d_{1}\int_{g_{\infty}}^{h_{\infty}}
J(x-y)\bar{w}(t,y)dy-d_{1}\bar{w}+\bar{w}(a_{T}-M_{1}w),\quad t>0,
x\in [g_{\infty}, h_{\infty}],\\[5pt]
\bar{w}(0,x)=u_{0}(x), ~x\in [-h_{0}, h_{0}];
\quad \bar{w}(0,x)=0,
~x\in [g_{\infty}, h_{\infty}]\setminus[-h_{0}, h_{0}].
\end{array}\right.
\end{align*}
Recall that in $(\ref{3.7})$ we have $\lambda_{1}(-(\mathcal{L}_{(g_{\infty},h_{\infty})}+a_{T}))
=\lambda_{1}(-(L_{(g_{\infty},h_{\infty})}+a))=0$, where $-(\mathcal{L}_{(g_{\infty},h_{\infty})}+a_{T})$
is defined in $(\ref{3.6})$. By Proposition 3.5 in \cite{cdll19} (see also \cite{baz07,cov10}), we know that
$\lim_{t\rightarrow \infty}\bar{w}(t,x)=0$ uniformly
in $[g_{\infty},h_{\infty}]$. Thus, $w(t,x)$ and
$\bar{u}(t,x)=e^{\int_{0}^{t}[a(s)-a_{T}]ds}w(t,x)$ converge to $0$ uniformly in
$[g_{\infty},h_{\infty}]$ as $t\rightarrow +\infty$, which implies that
$\lim_{t\rightarrow \infty}\bar{u}(t,x)=0$ uniformly in $[g_{\infty},h_{\infty}]$.

On the other hand, it is easy to know that
\begin{align*}
\left\{\begin{array}{l}
\bar{u}_t\geq d_{1}\int_{g(t)}^{h(t)}
J(x-y)\bar{u}(t,y)dy-d_{1}\bar{u}+\bar{u}(a(t)-\bar{u}),\quad t>0,
x\in (g(t), h(t)),\\[5pt]
\bar{u}(t,g(t))\geq 0,~ \bar{u}(t,h(t))\geq 0,\\[5pt]
\bar{u}(0,x)=u_{0}(x),~ x\in [-h_{0},h_{0}].
\end{array}\right.
\end{align*}
By the comparison principle (Lemma 2.2 in \cite{cdll19}), we know
$u(t,x)\leq \bar{u}(t,x)$ for any $t>0$ and $x\in [g(t), h(t)]$.
Thus, $\lim_{t\rightarrow \infty}\|u(t,\cdot)\|_{C([g(t),h(t)])}=0$.
\end{proof}

From Theorem $\ref{t4.4}$, we can obtain the following spreading-vanishing dichotomy.

\begin{corollary}\label{c4.5} (Spreading-vanishing dichotomy) Let $(u,v,g,h)$ be the unique solution
of $(\ref{1.1})$. Then, the following alternative holds:\\
\emph{Either} $(i)$ spreading: $\lim_{t\rightarrow \infty}(h(t)-g(t))=\infty$,
\emph{or} $(ii)$ vanishing: $\lim_{t\rightarrow \infty}(g(t), h(t))=(g_{\infty}, h_{\infty})$ is a finite
interval and $\lim_{t\rightarrow \infty}\max_{g(t)\leq x\leq h(t)}u(t,x)=\lim_{t\rightarrow \infty}\max_{g(t)\leq x\leq h(t)}\\v(t,x)=0$.
\end{corollary}

If we further assume the following weak competition condition
\begin{align*}\label{4.5}
\min_{[0,T]}a(t)>\max_{[0,T]}b(t)\cdot\max_{[0,T]}c(t), \quad
\min_{[0,T]}c(t)>\max_{[0,T]}d(t)\cdot\max_{[0,T]}a(t),
\tag{4.5}
\end{align*}
then we can establish the asymptotic estimates of $(u,v)$ when spreading occurs.
To achieve it, we first give a lemma concerning the asymptotic stability of
time-periodic solutions for the equations with nonlocal and mixed diffusions in $\mathbb{R}$.

\begin{lemma}\label{l4.6} $(i)$ For any bounded, uniformly continuous initial value $w_{0}$ with $\inf_{x\in \mathbb{R}}w_{0}>0$, the unique solution $w(t,x;w_{0})$ of
\begin{align*}
\left\{\begin{array}{l}
\partial_{t}w=d_{1}\left(\int_{\mathbb{R}}J(x-y)w(t,y)dy-w\right)+w(a(t)-w), \quad t>0,~x\in \mathbb{R},\\[5pt]
w(0,x)=w_{0}(x),\quad x\in \mathbb{R}
\end{array}\right.
\end{align*}
satisfies $\|w(t,\cdot;w_{0})-\phi^{*}(t)\|_{L^{\infty}(\mathbb{R})}\rightarrow 0$ as
$t\rightarrow \infty$, where $\phi^{*}(t)$ is the unique positive $T$-periodic solution of
\begin{align*}
\phi^{\prime}=\phi(a(t)-\phi), \quad
\phi(0)=\phi(T).
\end{align*}
$(ii)$ For any bounded, uniformly continuous initial value $w_{0}\in C^{2}(\mathbb{R})$ with $\inf_{x\in \mathbb{R}}w_{0}>0$, the unique solution $w(t,x;w_{0})$ of
\begin{align*}\label{4.6}
\left\{\begin{array}{l}
\partial_{t}w=d_{2}\left[\tau \partial_{x}^{2}w+(1-\tau)\left(\int_{\mathbb{R}}J(x-y)w(t,y)dy-w\right)\right]\\[3pt]
\qquad\quad+w(c(t)-w), \quad t>0,~x\in \mathbb{R},\\[3pt]
w(0,x)=w_{0}(x),\quad x\in \mathbb{R}
\end{array}\right.
\tag{4.6}
\end{align*}
satisfies $\|w(t,\cdot;w_{0})-\psi^{*}(t)\|_{L^{\infty}(\mathbb{R})}\rightarrow 0$ as
$t\rightarrow \infty$, where $\psi^{*}(t)$ is the unique positive $T$-periodic solution of
\begin{align*}
\psi^{\prime}=\psi(c(t)-\psi), \quad
\psi(0)=\psi(T).
\end{align*}
$(iii)$ The $T$-periodic functions $\phi^{*}(t)$, $\psi^{*}(t)$ in $(i)$-$(ii)$ are the unique positive solutions of
\begin{align*}\label{4.7}
\left\{\begin{array}{l}
\partial_{t}\phi=d_{1}\left(\int_{\mathbb{R}}J(x-y)\phi(t,y)dy-\phi\right)+\phi(a(t)-\phi), \quad t>0,~x\in \mathbb{R},\\[5pt]
\phi(0,x)=\phi(T,x),\quad x\in \mathbb{R}
\end{array}\right.
\tag{4.7}
\end{align*}
and
\begin{align*}\label{4.8}
\left\{\begin{array}{l}
\partial_{t}\psi=d_{2}\left[\tau \partial_{x}^{2}\psi+(1-\tau)\left(\int_{\mathbb{R}}J(x-y)\psi(t,y)dy-\psi\right)\right]\\[3pt]
\qquad\quad+\psi(c(t)-\psi), \quad t>0,~x\in \mathbb{R},\\[3pt]
\psi(0,x)=\psi(T,x),\quad x\in \mathbb{R},
\end{array}\right.
\tag{4.8}
\end{align*}
respectively.
\end{lemma}

\begin{proof} For $(i)$, the proof can be seen in that of Theorem 2.3(3) in \cite{ggs19}. We can apply similar arguments to prove $(ii)$, here we provide the details for the readers' convenient.

The existence and uniqueness of positive solution for $(\ref{4.6})$ can be established by applying
Lemma 2.3 in \cite{ch02} (the comparison principle) and the upper-lower solutions method.
Now we prove that $\psi^{*}(t)$ is globally asymptotically stable. By Example 4.19 in \cite{hk91}, we know that for any $\varepsilon>0$ and $K>0$, there exists $t_{\varepsilon, K}>0$ such that
\begin{align*}\label{4.9}
\Big|\psi^{*}(t)-\psi(t;K)\Big|\leq \varepsilon \quad \mbox{for}~t\geq t_{\varepsilon,K},
\tag{4.9}
\end{align*}
where $\psi(t;K)$ is the solution of
\begin{align*}
\psi^{\prime}=\psi(c(t)-\psi), \quad
\psi(0)=K.
\end{align*}

Denote $\underline{w}(0)=\inf_{x\in\mathbb{R}}w_{0}(x)$ and $\bar{w}(0)=\sup_{x\in\mathbb{R}}w_{0}(x)$.
Obviously, $\psi(t;\underline{w}(0))\leq \psi(t;\bar{w}(0))$.
From $(\ref{4.9})$, for any small $\varepsilon>0$, there exists $t_{\varepsilon}=\max\{t_{\varepsilon, \underline{w}(0)},t_{\varepsilon, \bar{w}(0)}\}>0$ such that
\begin{align*}\label{4.10}
0<\psi(t;\underline{w}(0))-\varepsilon\leq\psi^{*}(t)\leq\psi(t;\bar{w}(0))+\varepsilon.
\tag{4.10}
\end{align*}
Since
$(\psi(t;\underline{w}(0)), \psi(t;\bar{w}(0)))$ is a pair of upper-lower solutions of $(\ref{4.6})$, by Lemma 2.3 in \cite{ch02} (the comparison principle), we have
\begin{align*}\label{4.11}
\psi(t;\underline{w}(0))
\leq w(t,x;w_{0})
\leq \psi(t;\bar{w}(0)).
\tag{4.11}
\end{align*}

Next, we prove
\begin{align*}\label{4.12}
\begin{array}{rl}
0\leq \ln\frac{\psi(t;\bar{w}(0))}{\psi(t;\underline{w}(0))}\rightarrow 0
\quad \mbox{as}~t\rightarrow \infty.
\end{array}
\tag{4.12}
\end{align*}
In fact,
\begin{align*}
\begin{array}{rl}
\frac{d}{dt}\ln\frac{\psi(t;\bar{w}(0))}{\psi(t;\underline{w}(0))}
&=\frac{\psi^{\prime}(t;\bar{w}(0))}{\psi(t;\bar{w}(0))}
-\frac{\psi^{\prime}(t;\underline{w}(0))}{\psi(t;\underline{w}(0))}
=-(\psi(t;\bar{w}(0))-\psi(t;\underline{w}(0)))\\[3pt]
&\leq -\psi(t;\underline{w}(0))\ln\frac{\psi(t;\bar{w}(0))}{\psi(t;\underline{w}(0))},
\end{array}
\end{align*}
where we have used the inequality $-(a-b)\leq -b\ln\frac{a}{b}$ for $a>b>0$.
Let $K_{0}:=\inf_{t\geq 0}\psi(t;\underline{w}(0))>0$. We have
\begin{align*}
\begin{array}{rl}
\frac{d}{dt}\ln\frac{\psi(t;\bar{w}(0))}{\psi(t;\underline{w}(0))}
\leq -K_{0}\ln\frac{\psi(t;\bar{w}(0))}{\psi(t;\underline{w}(0))},
\end{array}
\end{align*}
which implies
\begin{align*}
\begin{array}{rl}
0\leq \ln\frac{\psi(t;\bar{w}(0))}{\psi(t;\underline{w}(0))}
\leq \ln\frac{\bar{w}(0)}{\underline{w}(0)}e^{-K_{0}t}\rightarrow 0
\quad \mbox{as}~t\rightarrow \infty.
\end{array}
\end{align*}
Combing $(\ref{4.10})$-$(\ref{4.12})$, we get the result of $(ii)$.

For the uniqueness of solution to $(\ref{4.7})$ in $(iii)$, we first give a lower bound estimate of any bounded positive solutions of $(\ref{4.7})$. Consider the following stationary problem
\begin{align*}\label{4.13}
\begin{array}{rl}
-d_{1}\left(\int_{-l}^{l}J(x-y)\phi(y)dy-\phi(x)\right)=\phi(\min_{t\in[0,T]} a(t)-\phi), \quad -l<x<l.
\tag{4.13}
\end{array}
\end{align*}
From Proposition 3.6 in \cite{cdll19}, $(\ref{4.13})$ has a unique positive bounded solution $\phi^{l}(x)$ for sufficiently large $l$, and
$\phi^{l}(x)\rightarrow \min_{t\in[0,T]} a(t)$ in $C([-L,L])$ for any $L>0$ as
$l\rightarrow +\infty$.
For any positive solution $\hat{\phi}(t,x)$ of (4.7), by the comparison principle (Lemma $\ref{l3-3}$), we have $\hat{\phi}\geq \phi^{l}$ on $[0,T]\times [-l,l]$.
Letting $l\rightarrow +\infty$, we get $\hat{\phi}\geq \min_{t\in[0,T]} a(t)>0$.

Next, we briefly show that $(\ref{4.7})$ has a minimal positive solution. Consider the problem
\begin{align*}\label{414}
\left\{\begin{array}{l}
\partial_{t}\phi=d_{1}\left(\int_{-l}^{l}J(x-y)\phi(t,y)dy-\phi\right)+\phi(a(t)-\phi), \quad 0\leq t\leq T,~-l<x<l,\\[5pt]
\phi(0,x)=\phi(T,x),\quad -l\leq x\leq l.
\end{array}\right.
\tag{4.14}
\end{align*}
For sufficiently large $l$, $(\ref{414})$ admits a unique positive solution $\psi_{*}^{l}$. By the comparison principle (Lemma $\ref{l3-3}$), we can show that $\psi_{*}^{l}$ is increasing in $l$ and $\psi_{*}^{l}\leq \hat{\phi}(t,x)$ on $[0,T]\times [-l,l]$ for any positive solution $\hat{\phi}(t,x)$ of $(\ref{4.7})$ and any $l>0$.
Thus, the limit function $\psi_{*}=\lim_{l\rightarrow \infty}\psi_{*}^{l}$ is exactly a minimal positive solution of $(\ref{4.7})$.

Finally, we prove the uniqueness by using a technique introduced by Marcus \& V\'{e}ron \cite{mv98}.
Arguing indirectly, we assume that $(\ref{4.7})$ has a positive bounded solution $\hat{\psi}$
such that $\hat{\psi}\not\equiv \psi_{*}$. Then there exists a constant $k>1$ such that
$\psi_{*}\leq \hat{\psi}\leq k\psi_{*}$ in $[0,T]\times \mathbb{R}$. By the strong maximum principle
(see Definition 1.4 and Theorem \textbf{F} in \cite{svo19}), we have $\psi_{*}<\hat{\psi}$. Define
$\bar{\psi}=\psi_{*}-(2k)^{-1}(\hat{\psi}-\psi_{*})$. By direct calculations,
we get
\begin{align*}\label{4.15}
\psi_{*}>\bar{\psi}\geq\frac{k+1}{2k}\psi_{*},\quad
\frac{2k}{2k+1}\bar{\psi}+\frac{1}{2k+1}\hat{\psi}=\psi_{*}.
\tag{4.15}
\end{align*}
Note that $\hat{\psi}(\hat{\psi}-a(t))$ is convex in $\hat{\psi}\in (0,+\infty)$.
We have $\psi_{*}^{2}\leq \frac{2k}{2k+1}\bar{\psi}^{2}+\frac{1}{2k+1}\hat{\psi}^{2}$.
It is easy to check that
\begin{align*}
\begin{array}{rl}
\partial_{t}\bar{\psi}
&\geq d_{1}\left(\int_{\mathbb{R}}J(x-y)\bar{\psi}(t,y)dy-\bar{\psi}\right)+\bar{\psi}(a(t)-\bar{\psi})\\[3pt]
&\geq d_{1}\left(\int_{-l}^{l}J(x-y)\bar{\psi}(t,y)dy-\bar{\psi}\right)+\bar{\psi}(a(t)-\bar{\psi}), \quad 0\leq t\leq T,~x\in [-l,l]
\end{array}
\end{align*}
and $\bar{\psi}(0,x)=\bar{\psi}(T,x)$ for $x\in[-l,l]$.
Thus, $\bar{\psi}$ is an upper solution of $(\ref{414})$. By the comparison principle (Lemma $\ref{l3-3}$), we have
$\psi_{*}^{l}\leq \bar{\psi}$ in $[0,T]\times [-l,l]$. Since $\psi_{*}^{l}\rightarrow \psi_{*}$
in $C^{1,0}([0,T]\times[-L,L])$ for any $L>0$ as $l\rightarrow\infty$. It follows that
$\psi_{*}\leq \bar{\psi}$ in $[0,T]\times\mathbb{R}$,
which contradicts with $(\ref{4.15})$. This completes the proof of the uniqueness of positive solution to $(\ref{4.7})$.

In proving the uniqueness of positive solution to $(\ref{4.8})$, we need to replace the auxiliary problems $(\ref{4.13})$ and $(\ref{414})$ with
\begin{align*}
\left\{\begin{array}{l}
-d_{2}\left[\tau\partial_{x}^{2}\phi+(1-\tau)\left(\int_{-l}^{l}J(x-y)\phi(y)dy-\phi(x)\right)\right]=\phi(\min_{t\in[0,T]} c(t)-\phi), \quad -l<x<l,\\[3pt]
\phi(\pm l)=0
\end{array}\right.
\end{align*}
and
\begin{align*}
\left\{\begin{array}{l}
\partial_{t}\phi=d_{2}\left[\tau\partial_{x}^{2}\phi+(1-\tau)\left(\int_{-l}^{l}J(x-y)\phi(t,y)dy-\phi(t,x)\right)\right]\\[3pt]
\qquad\quad+\phi(c(t)-\phi), \quad 0\leq t\leq T,~-l<x<l,\\[3pt]
\phi(t,\pm l)=0,\quad 0\leq t\leq T,\\[3pt]
\phi(0,x)=\phi(T,x),\quad -l\leq x\leq l,
\end{array}\right.
\end{align*}
respectively, which have zero boundary conditions. The proof is similar as that of $(\ref{4.7})$, here we omit the details.
\end{proof}

\begin{theorem}\label{t4.7} Suppose that $(\ref{4.5})$ holds and $h_{\infty}-g_{\infty}=\infty$.
Then
\begin{align*}
\begin{array}{rl}
&U_{*}(t)\leq\liminf_{n\rightarrow \infty}u(t+nT, x)\leq\limsup_{n\rightarrow \infty}u(t+nT, x)\leq U^{*}(t),\\[3pt]
&V_{*}(t)\leq\liminf_{n\rightarrow \infty}v(t+nT, x)\leq\limsup_{n\rightarrow \infty}v(t+nT, x)\leq V^{*}(t)
\end{array}
\end{align*}
uniformly in $[0,T]\times[-L,L]$ for any $L>0$, where $U^{*}(t), V^{*}(t)$,
$V_{*}(t)$ and $U_{*}(t)$ are positive $T$-periodic solutions of
\begin{align*}
\begin{array}{rl}
&\frac{d U^{*}}{dt}=U^{*}(a(t)-U^{*}),\quad U^{*}(0)=U^{*}(T),\\[3pt]
&\frac{d V^{*}}{dt}=V^{*}(c(t)-V^{*}),\quad U^{*}(0)=U^{*}(T),\\[3pt]
&\frac{d V_{*}}{dt}=V_{*}\Big(c(t)-V_{*}-d(t)U^{*}(t)\Big),\quad V_{*}(0)=V_{*}(T)
\end{array}
\end{align*}
and
\begin{align*}
\begin{array}{rl}
\frac{d U_{*}}{dt}=U_{*}\Big(a(t)-U_{*}-b(t)V^{*}(t)\Big),\quad V_{*}(0)=V_{*}(T),
\end{array}
\end{align*}
respectively.
\end{theorem}

\begin{proof} In Theorem 3.2 of \cite{waz16}, similar results have been obtained for the random dispersal case. Since the nonlocal dispersal is considered here, we give the details.

\emph{Step 1.} $\limsup_{n\rightarrow \infty}u(t+nT,x)\leq U^{*}(t)$, $\limsup_{n\rightarrow \infty}v(t+nT,x)\leq V^{*}(t)$ uniformly in $[0,T]\times [-L,L]$ for any given $L>0$.

Let $w(t,x)$ be the unique positive solution of
\begin{align*}
\left\{\begin{array}{l}
\partial_{t}w=d_{1}\left(\int_{\mathbb{R}}J(x-y)w(t,y)dy-w\right)+w(a(t)-w), \quad t>0,~x\in \mathbb{R},\\[5pt]
w(0,x)=\|u_{0}\|_{L^{\infty}([-h_{0},h_{0}])}>0,\quad x\in \mathbb{R}.
\end{array}\right.
\end{align*}
By Lemma $\ref{l4.6}$$(i)$, we know that $\lim_{n\rightarrow \infty}w(t+nT,x)\rightarrow U^{*}(t)$ uniformly
for $(t,x)\in [0,T]\times [-L,L]$. Moreover, since $w$ satisfies
\begin{align*}
\begin{array}{rl}
\partial_{t}w\geq d_{1}\left(\int_{g(t)}^{h(t)}J(x-y)w(t,y)dy-w\right)+w(a(t)-w), \quad t>0,~x\in (g(t),h(t)),
\end{array}
\end{align*}
by the comparison principle (Lemma 2.2 in \cite{cdll19})
we have $u(t,x)\leq w(t,x)$ for $(t,x)\in[0,+\infty)\times [g(t),h(t)]$.
Thus, $\limsup_{n\rightarrow \infty}u(t+nT,x)\leq U^{*}(t)$ uniformly in $[0,T]\times [-L,L]$.

Similarly, by applying Lemma $\ref{l4.6}$$(ii)$ and Lemma $\ref{l2.2}$, we can prove $\limsup_{n\rightarrow \infty}v(t+nT,x)\leq V^{*}(t)$ uniformly in $[0,T]\times [-L,L]$.

\emph{Step 2.} $\liminf_{n\rightarrow \infty}v(t+nT,x)\geq V_{*}(t)$ uniformly in $[0,T]\times [-L,L]$ for any given $L>0$.

By the assumption $(\ref{4.5})$ and the fact $U^{*}\leq \max_{[0,T]}a(t)$, we know that there exists $\varepsilon_{0}>0$ such that
\begin{align*}
c_{\varepsilon}(t):=c(t)-d(t)(U^{*}(t)+\varepsilon)\geq\min_{[0,T]}c(t)-\max_{[0,T]}d(t)\cdot(\max_{[0,T]}a(t)+\varepsilon)>0 \end{align*}
for any $0<\varepsilon\leq \varepsilon_{0}$. For such a fixed $\varepsilon$, from Lemma $\ref{l3.7}$$(i)$ we can deduce that there exists $L_{\varepsilon}>L$ such that $\lambda_{1}(-(\tilde{L}_{(-l, l)}+c_{\varepsilon}))<0$ for all $l\geq L_{\varepsilon}$. Since
$h_{\infty}-g_{\infty}=\infty$ and $\limsup_{n\rightarrow \infty}u(t+nT,x)\leq U^{*}(t)$ locally uniformly in $[0,T]\times \mathbb{R}$, for any fixed $\varepsilon\in (0,\varepsilon_{0})$ and $l>L_{\varepsilon}$ there exists $m\in \mathbb{N}$ such that
\begin{align*}
g(t)<-l,~h(t)>l,~u(t,x)<U^{*}(t)+\varepsilon,\quad \forall t\geq mT,~-l\leq x\leq l.
\end{align*}

Let $z_{l}^{\varepsilon}$ be the unique positive solution of
\begin{align*}
\left\{\begin{array}{l}
\partial_{t}z=d_{2}\left[\tau \partial_{x}^{2}z+(1-\tau)\left(\int_{-l}^{l}J(x-y)z(t,y)dy-z\right)\right]\\[3pt]
\qquad\quad+z(c_{\varepsilon}(t)-z), \quad t>mT,~-l<x<l,\\[3pt]
z(t,\pm l)=0, \quad t>mT,\\[3pt]
z(mT,x)=v(mT,x), \quad -l<x<l.
\end{array}\right.
\end{align*}
By the comparison principle derived from Lemma $\ref{l2.2}$, $v(t,x)\geq z_{l}^{\varepsilon}(t,x)$ for
$t\geq mT$ and $x\in[-l,l]$. Since $\lambda_{1}(-(\tilde{L}_{(-l,l)}+c_{\varepsilon}))<0$, by Lemma $\ref{l3.7}$$(ii)$, we deduce that $\lim_{n\rightarrow \infty}z_{l}^{\varepsilon}(t+nT,x)
=Z_{l}^{\varepsilon}(t,x)$ in $C^{1,2}([0,T]\times [-l,l])$, where $Z_{l}^{\varepsilon}(t,x)$
is the unique positive $T$-periodic solution of
\begin{align*}
\left\{\begin{array}{l}
\partial_{t}Z=d_{2}\left[\tau \partial_{x}^{2}Z+(1-\tau)\left(\int_{-l}^{l}J(x-y)Z(t,y)dy-Z\right)\right]\\[3pt]
\qquad\quad+Z(c_{\varepsilon}(t)-Z), \quad 0\leq t\leq T,~-l<x<l,\\[3pt]
Z(t,\pm l)=0, \quad 0\leq t\leq T,\\[3pt]
Z(0,x)=Z(T,x), \quad -l<x<l.
\end{array}\right.
\end{align*}
By the comparison principle, we can prove that $Z_{l}^{\varepsilon}(t,x)$ is increasing
with respect to $l$. Thus,
\begin{align*}
\lim_{l\rightarrow +\infty}Z_{l}^{\varepsilon}(t,x)=Z^{\varepsilon}(t,x)\quad
\mbox{in}~C^{1,2}([0,T]\times[-L,L]),
\end{align*}
where $Z^{\varepsilon}(t,x)$ is the unique positive $T$-periodic solution of
\begin{align*}
\left\{\begin{array}{l}
\partial_{t}Z=d_{2}\left[\tau \partial_{x}^{2}Z+(1-\tau)\left(\int_{\mathbb{R}}J(x-y)Z(t,y)dy-Z\right)\right]\\[3pt]
\qquad\quad+Z(c_{\varepsilon}(t)-Z), \quad t\in[0,T],~x\in \mathbb{R},\\[3pt]
Z(0,x)=Z(T,x), \quad x\in \mathbb{R}.
\end{array}\right.
\end{align*}
By Lemma $\ref{l4.6}$$(iii)$, $Z^{\varepsilon}$ satisfies
\begin{align*}
\frac{dZ}{dt}=Z(c_{\varepsilon}(t)-Z), \quad
Z(0)=Z(T).
\end{align*}
Thus, we have
$\lim_{n\rightarrow \infty}v(t+nT,x)
\geq Z^{\varepsilon}(t,x)$ uniformly for $(t,x)\in[0,T]\times[-L,L]$. Letting
$\varepsilon\rightarrow 0$, we know that $\liminf_{n\rightarrow \infty}v(t+nT,x)\geq V_{*}(t)$ uniformly in $[0,T]\times [-L,L]$ for any given $L>0$.

\emph{Step 3.} $\liminf_{n\rightarrow \infty}u(t+nT,x)\geq U_{*}(t)$ uniformly in $[0,T]\times [-L,L]$ for any given $L>0$.

By the assumption $(\ref{4.5})$ and the fact $V^{*}\leq \max_{[0,T]}c(t)$, we know that there exists $\varepsilon_{1}>0$ such that
\begin{align*}
a_{\varepsilon}(t):=a(t)-b(t)(V^{*}(t)+\varepsilon)\geq\min_{[0,T]}a(t)-\max_{[0,T]}b(t)\cdot(\max_{[0,T]}c(t)+\varepsilon)>0 \end{align*}
for any $0<\varepsilon\leq \varepsilon_{1}$. For such a fixed $\varepsilon$, from Lemma $\ref{l3.6}$ we can deduce that there exists $l_{\varepsilon}>L$ such that $\lambda_{1}(-(L_{(-l,l)}+a_{\varepsilon}))<0$ for all $l\geq l_{\varepsilon}$. Since
$h_{\infty}-g_{\infty}=\infty$ and $\limsup_{n\rightarrow \infty}v(t+nT,x)\leq V^{*}(t)$ locally uniformly in $[0,T]\times \mathbb{R}$, for any fixed $\varepsilon\in (0,\varepsilon_{1})$ and $l>l_{\varepsilon}$ there exists $m_{1}\in \mathbb{N}$ such that
\begin{align*}
g(t)<-l,~h(t)>l,~v(t,x)<V^{*}(t)+\varepsilon,\quad \forall t\geq m_{1}T,~-l\leq x\leq l.
\end{align*}

Let $p_{l}^{\varepsilon}$ be the unique positive solution of
\begin{align*}
\left\{\begin{array}{l}
\partial_{t}p=d_{1}\left(\int_{-l}^{l}J(x-y)p(t,y)dy-p\right)
+p(a_{\varepsilon}(t)-p), \quad t>m_{1}T,~-l<x<l,\\[5pt]
p(m_{1}T,x)=u(m_{1}T,x), \quad -l<x<l.
\end{array}\right.
\end{align*}
Since $\lambda_{1}(-(L_{(-l,l)}+a_{\varepsilon}))<0$, by Lemma $\ref{l3.4}$, we know that $\lim_{n\rightarrow \infty}p_{l}^{\varepsilon}(t+nT,x)
=P_{l}^{\varepsilon}(t,x)$ in $C^{1,0}([0,T]\times [-l,l])$, where $P_{l}^{\varepsilon}(t,x)$
is the unique positive $T$-periodic solution of
\begin{align*}
\left\{\begin{array}{l}
\partial_{t}P=d_{1}\left(\int_{-l}^{l}J(x-y)P(t,y)dy-P\right)
+P(a_{\varepsilon}(t)-P), \quad 0\leq t\leq T,~-l<x<l,\\[5pt]
P(0,x)=P(T,x), \quad -l<x<l.
\end{array}\right.
\end{align*}
By the comparison principle (Lemma $\ref{l3-3}$), $P_{l}^{\varepsilon}(t,x)$ is increasing in $l$. Thus,
\begin{align*}
\lim_{l\rightarrow +\infty}P_{l}^{\varepsilon}(t,x)=P^{\varepsilon}(t,x)\quad
\mbox{in}~C^{1,0}([0,T]\times[-L,L]),
\end{align*}
where $P^{\varepsilon}(t,x)$ is the unique positive $T$-periodic solution of
\begin{align*}
\left\{\begin{array}{l}
\partial_{t}P=d_{1}\left(\int_{\mathbb{R}}J(x-y)P(t,y)dy-P\right)
+P(a_{\varepsilon}(t)-P), \quad t\in[0,T],~x\in \mathbb{R},\\[5pt]
P(0,x)=P(T,x), \quad x\in \mathbb{R}.
\end{array}\right.
\end{align*}
By Lemma $\ref{l4.6}$$(iii)$, $P^{\varepsilon}$ satisfies
\begin{align*}
\frac{dP}{dt}=P(a_{\varepsilon}(t)-P), \quad
P(0)=P(T).
\end{align*}
Thus, we have
$\lim_{n\rightarrow \infty}u(t+nT,x)
\geq P^{\varepsilon}(t,x)$ uniformly for $(t,x)\in[0,T]\times[-L,L]$. Letting
$\varepsilon\rightarrow 0$, we know that $\liminf_{n\rightarrow \infty}u(t+nT,x)\geq U_{*}(t)$ uniformly in $[0,T]\times [-L,L]$ for any given $L>0$.
\end{proof}

In what follows, we will provide some sufficient conditions for spreading and vanishing.

\begin{theorem}\label{t4.8} If $h_{\infty}-g_{\infty}<\infty$,
then
$h_{\infty}-g_{\infty}\leq h^{*}$,
where $|\Omega|=h^{*}$ is the unique root of $\lambda_{1}(-(\tilde{L}_{\Omega}+c))=0$
with $-(\tilde{L}_{\Omega}+c)$ defined as in $(\ref{3.8})$.
\end{theorem}

\begin{proof} Recall that in Theorem $\ref{t4.4}$ we have showed
that $h_{\infty}-g_{\infty}<\infty$ implies
\begin{align*}\label{4.16}
\lim_{t\rightarrow \infty}\|u(t,\cdot)\|_{C([g(t),h(t)])}=
\lim_{t\rightarrow \infty}\|v(t,\cdot)\|_{C([g(t),h(t)])}=0.
\tag{4.16}
\end{align*}
Assume on the contrary that $h_{\infty}-g_{\infty}>h^{*}$.
Then there exists $0<\varepsilon\ll 1$ and $\mathcal{T}\gg 1$ such that
\begin{align*}
\begin{array}{l}
h_{\infty}^{-\varepsilon}-g_{\infty}^{+\varepsilon}
=h_{\infty}-g_{\infty}-2\varepsilon>h^{*}_{\varepsilon},\\[3pt]
g(\mathcal{T})<g_{\infty}^{+\varepsilon},~h(\mathcal{T})>h_{\infty}^{-\varepsilon},\\[3pt]
0\leq u(t,x)<\varepsilon,~\forall t\geq \mathcal{T}, x\in [g_{\infty}^{+\varepsilon},h_{\infty}^{-\varepsilon}],
\end{array}
\end{align*}
where $|\Omega|=h^{*}_{\varepsilon}$ is the unique root of $\lambda_{1}(-(\tilde{L}_{\Omega}+c(t)-d(t)\varepsilon))=0$.
Then $v$ satisfies
\begin{align*}
\left\{\begin{array}{l}
v_t\geq d_{2}\left[\tau v_{xx}+(1-\tau)\left(\int_{g_{\infty}^{+\varepsilon}}^{h_{\infty}^{-\varepsilon}}J(x-y)v(t,y)dy-v\right)\right]\\[3pt]
\qquad+v(c(t)-d(t)\varepsilon-v),
\quad t>\mathcal{T},~
x\in (g_{\infty}^{+\varepsilon}, h_{\infty}^{-\varepsilon}),\\[3pt]
v(t,g_{\infty}^{+\varepsilon})> 0,~ v(t,h_{\infty}^{-\varepsilon})> 0, ~ t\geq\mathcal{T},\\[3pt]
v(\mathcal{T},x)>0,~ x\in (g_{\infty}^{+\varepsilon},h_{\infty}^{-\varepsilon}).
\end{array}\right.
\end{align*}
Let $\psi$ be the unique positive solution of
\begin{align*}
\left\{\begin{array}{l}
\psi_t=d_{2}\left[\tau \psi_{xx}+(1-\tau)\left(\int_{g_{\infty}^{+\varepsilon}}^{h_{\infty}^{-\varepsilon}}J(x-y)\psi(t,y)dy-\psi\right)\right]\\[3pt]
\qquad+\psi(c(t)-d(t)\varepsilon-\psi),
\quad
t>\mathcal{T},~
x\in (g_{\infty}^{+\varepsilon}, h_{\infty}^{-\varepsilon}),\\[3pt]
\psi(t,g_{\infty}^{+\varepsilon})=0,~ \psi(t,h_{\infty}^{-\varepsilon})=0, ~ t\geq\mathcal{T},\\[3pt]
\psi(\mathcal{T},x)=v(\mathcal{T},x),~ x\in (g_{\infty}^{+\varepsilon},h_{\infty}^{-\varepsilon}).
\end{array}\right.
\end{align*}
By Lemma $\ref{l2.2}$, we have
\begin{align*}
\psi(t,x)\leq v(t,x),\quad
t\geq \mathcal{T}, x\in [g_{\infty}^{+\varepsilon},h_{\infty}^{-\varepsilon}].
\end{align*}
Since $h_{\infty}^{-\varepsilon}-g_{\infty}^{+\varepsilon}=h_{\infty}-g_{\infty}-2\varepsilon>h^{*}_{\varepsilon}$, we have
$\lambda_{1}(-(\tilde{L}_{(g_{\infty}^{+\varepsilon}, h_{\infty}^{-\varepsilon})}+c(t)-d(t)\varepsilon))<0$, and then Lemma $\ref{l3.7}$$(ii)$ implies that
$\psi(t+nT,x)\rightarrow \omega(t,x)$ as $n\rightarrow\infty$ uniformly in
the compact subset of $(g_{\infty}^{+\varepsilon},h_{\infty}^{-\varepsilon})$,
where $\omega(t,x)$ is the unique positive periodic solution of
\begin{align*}
\left\{\begin{array}{l}
\omega_t=d_{2}\left[\tau \omega_{xx}+(1-\tau)\left(\int_{g_{\infty}^{+\varepsilon}}^{h_{\infty}^{-\varepsilon}}J(x-y)\omega(t,y)dy-\omega\right)\right]\\[3pt]
\qquad+\omega(c(t)-d(t)\varepsilon-\omega),\quad
t\in [0, T],~x\in (g_{\infty}^{+\varepsilon}, h_{\infty}^{-\varepsilon}),\\[3pt]
\omega(t,g_{\infty}^{+\varepsilon})=0,~ \omega(t,h_{\infty}^{-\varepsilon})=0, \quad t\in [0, T],\\[3pt]
\omega(0,x)=\omega(T,x), \quad x\in (g_{\infty}^{+\varepsilon},h_{\infty}^{-\varepsilon}).
\end{array}\right.
\end{align*}
Therefore, $\liminf_{n\rightarrow \infty}v(t+nT,x)\geq \lim_{n\rightarrow \infty}\psi(t+nT,x)=\omega(t,x)>0$ for all $x\in (g_{\infty}^{+\varepsilon}, h_{\infty}^{-\varepsilon})$, which is a contradiction to $(\ref{4.16})$. This completes the proof.
\end{proof}

\begin{corollary}\label{r4.9} If $h_{0}\geq\frac{1}{2}h^{*}$,
then spreading occurs, that is, $h_{\infty}-g_{\infty}=+\infty$.
\end{corollary}

If $a_{T}\geq d_{1}$, then Lemma $\ref{l3.6}$ implies that $\lambda_{1}(-(L_{\Omega}+a))<0$ for all $l:=|\Omega|>0$. Thus, the vanishing can not happen by the proof of Theorem $\ref{t4.4}$, which means that
$h_{\infty}-g_{\infty}=+\infty$ always holds.

\begin{theorem}\label{t4.10}
If $a_{T}\geq d_{1}$, then spreading always happens.
\end{theorem}

On the other hand, if $a_{T}<d_{1}$, then Lemma $\ref{l3.6}$ implies that $\lambda_{1}(-(L_{\Omega}+a))>0$
for $0<|\Omega|\ll 1$, and $\lambda_{1}(-(L_{\Omega}+a))<0$
for $|\Omega|\gg 1$. Since $\lambda_{1}(-(L_{\Omega}+a))$ is strictly decreasing in $|\Omega|$, there exists a $l^{*}>0$ such that $\lambda_{1}(-(L_{\Omega}+a))=0$ for $|\Omega|=l^{*}$, $\lambda_{1}(-(L_{\Omega}+a))>0$
for $|\Omega|<l^{*}$ and $\lambda_{1}(-(L_{\Omega}+a))<0$ for $|\Omega|>l^{*}$. From the proof of $(\ref{4.3})$, we know that if $h_{\infty}-g_{\infty}<+\infty$ then $h_{\infty}-g_{\infty}\leq l^{*}$. Therefore, if
$h_{0}\geq \frac{l^{*}}{2}$ then we have $h_{\infty}-g_{\infty}=+\infty$.

\begin{theorem}\label{t4.11} Assume $a_{T}<d_{1}$ and $h_{0}<\frac{1}{2}\min\{h^{*}, l^{*}\}$.
Then there exists $\Lambda_{0}>0$ such that $h_{\infty}-g_{\infty}<+\infty$ when
$\mu+\rho_{1}+\rho_{2}\leq \Lambda_{0}$.
\end{theorem}

\begin{proof}
Since $\lambda_{1}(-(L_{(-h_{0}, h_{0})}+a))>0$, we can choose $h_{0}<h_{1}<\frac{l^{*}}{2}$
such that $\lambda:=\lambda_{1}(-(L_{(-h_{1}, h_{1})}+a))>0$.

Let $\bar{u}$ be the unique solution of
\begin{align*}
\left\{\begin{array}{l}
\bar{u}_t=d_{1}\int_{-h_{1}}^{h_{1}}
J(x-y)\bar{u}(t,y)dy-d_{1}\bar{u}+a(t)\bar{u},\quad t>0,
x\in [-h_{1}, h_{1}],\\[5pt]
\bar{u}(0,x)=u_{0}(x),\quad |x|\leq h_{0};\quad
\bar{u}(0,x)=0,\quad h_{0}<|x|\leq h_{1}.
\end{array}\right.
\end{align*}
And let $\varphi(t,x)$ be the corresponding eigenfunction associated with $\lambda$ and satisfies
$\|\varphi\|_{L^{\infty}([0, T]\times [-h_{1}, h_{1}])}=1$, that is,
\begin{align*}
-\left(L_{(-h_{1}, h_{1})}+a\right)[\varphi]
=\lambda \varphi.
\end{align*}
Let $\omega(t,x)=Ce^{-\frac{\lambda t}{2}}\varphi(t,x)$ for some $C>0$, it is easy to check that
\begin{align*}
\begin{array}{rl}
&\omega_{t}-d_{1}\int_{-h_{1}}^{h_{1}}
J(x-y)\omega(t,y)dy+d_{1}\omega-a(t)\omega\\[3pt]
&=Ce^{-\frac{\lambda t}{2}}\left(\varphi_{t}-d_{1}\int_{-h_{1}}^{h_{1}}
J(x-y)\varphi(t,y)dy+d_{1}\varphi-a(t)\varphi-\frac{\lambda}{2}\varphi\right)\\[3pt]
&=\frac{1}{2}\lambda Ce^{-\frac{\lambda t}{2}}\varphi(t,x)>0,
\end{array}
\end{align*}
for all $t>0$ and $x\in [-h_{1}, h_{1}]$. Choosing $C>0$ large such that $\omega(0,x)=C\varphi(0,x)>u_{0}(x)$ on $[-h_{1}, h_{1}]$. Applying Lemma 3.3 in \cite{cdll19},
we have
\begin{align*}
\bar{u}(t,x)\leq \omega(t,x)=Ce^{-\frac{\lambda t}{2}}\varphi(t,x)
\leq Ce^{-\frac{\lambda t}{2}},
\end{align*}
for all $t>0$ and $x\in [-h_{1}, h_{1}]$.

On the other hand, since $h_{0}<\frac{h^{*}}{2}$, we have $\lambda_{1}(-(\tilde{L}_{(-h_{0}, h_{0})}+c))>0$. Then there exists $0<\varepsilon_{0}\ll 1$ such that for any $0<\varepsilon\leq \varepsilon_{0}$, the following eigenvalue problem has a postive principle eigenvalue $\tilde{\lambda}_{1}>0$:
\begin{align*}
\left\{\begin{array}{l}
\varphi_{t}-d_{2}\Big[\tau\varphi_{xx}+(1-\tau)\Big((1+\varepsilon)\int_{\Omega}J(x-y)\varphi(t,y)dy-\varphi(t,x)\Big)\Big]-c(t)\varphi=\lambda \varphi \\[3pt]
\qquad\qquad\qquad\qquad\qquad\qquad\qquad\qquad\qquad\qquad\qquad\qquad\quad
\mbox{in}~[0,T]\times(-h_0, h_0),\\[3pt]
\varphi(t,\pm h_{0})=0 \quad \mbox{in}~[0,T], \\[3pt]
\varphi(0,x)=\varphi(T,x) \quad \mbox{in}~[-h_0,h_0].
\end{array}\right.
\end{align*}
Let $\tilde{\varphi}(t,x)$ be the corresponding normalized eigenfunction associated with $\tilde{\lambda}_{1}$. Note that $\tilde{\varphi}_{x}(t,h_{0})<0, \tilde{\varphi}_{x}(t,-h_{0})>0$ in $[0,T]$. Then
there exists a constant $\alpha>0$ such that
\begin{align*}
x\tilde{\varphi}_{x}(t,x)\leq \alpha\tilde{\varphi}(t,x), \quad \forall(t,x)\in [0,T]\times [-h_{0},h_{0}].
\end{align*}

For any $(t,x)\in [0,\infty)\times [-s(t),s(t)]$, we define
\begin{align*}
\begin{array}{rl}
s(t)=h_{0}\varsigma(t), \quad \varsigma(t)=1+2\delta-\delta e^{-\sigma t}, \quad
\bar{v}(t,x)=ke^{-\sigma t}\tilde{\varphi}(\xi(t), \eta(t,x))
\end{array}
\end{align*}
with
\begin{align*}
\begin{array}{rl}
\xi(t)=\int_{0}^{t}\frac{1}{\varsigma^{2}(\theta)}d\theta,
\quad \eta(t,x)=\frac{h_{0}}{s(t)}x=\frac{x}{\varsigma(t)},
\end{array}
\end{align*}
where $k>0, \sigma>0,0<\delta<\frac{1}{2}(\frac{h_{1}}{h_{0}}-1)$ are positive constants to be determined later.
Then $\overline{v}(t,x)$ satisfies
\begin{align*}
\begin{array}{rl}
&\bar{v}_{t}(t,x)-d_{2}[\tau\bar{v}_{xx}+(1-\tau)(\int_{-s(t)}^{s(t)}J(x-y)\bar{v}(t,y)dy-\bar{v}(t,x))]
-\bar{v}(t,x)(c(t)-\bar{v}(t,x))\\[3pt]
&=ke^{-\sigma t}\Big[-\sigma\tilde{\varphi}(\xi,\eta)
-\frac{\varsigma^{\prime}(t)}{\varsigma(t)}\eta\tilde{\varphi}_{\eta}(\xi,\eta)
+d_{2}(1-\tau)\Big(\frac{1+\varepsilon}{\varsigma^{2}(t)}\int_{-h_{0}}^{h_{0}}J(\eta-\tilde{\eta})\tilde{\varphi}(\xi,\tilde{\eta})d\tilde{\eta}\\[3pt]
&\qquad\qquad-\varsigma(t)\int_{-h_{0}}^{h_{0}}J(\varsigma(t)\eta-\varsigma(t)\tilde{\eta})\tilde{\varphi}(\xi,\tilde{\eta})d\tilde{\eta}\Big)
+d_{2}(1-\tau)(1-\frac{1}{\varsigma^{2}(t)})\tilde{\varphi}(\xi,\eta)\\[3pt]
&\qquad\qquad+(\frac{1}{\varsigma^{2}(t)}c(\xi)-c(t))\tilde{\varphi}(\xi,\eta)
+\frac{1}{\varsigma^{2}(t)}\tilde{\lambda}_{1}\tilde{\varphi}(\xi,\eta)+ke^{-\sigma t}\tilde{\varphi}^{2}(\xi,\eta)\Big]\\[3pt]
&\geq ke^{-\sigma t}\Big[\Big(-\sigma-\sigma \alpha
+d_{2}(1-\tau)(1-\frac{1}{\varsigma^{2}(t)})+\frac{1}{\varsigma^{2}(t)}\tilde{\lambda}_{1}
+(\frac{1}{\varsigma^{2}(t)}c(\xi)-c(t))\Big)\tilde{\varphi}(\xi,\eta)\\[3pt]
&\quad+d_{2}(1-\tau)\Big(\frac{1+\varepsilon}{\varsigma^{2}(t)}\int_{-h_{0}}^{h_{0}}J(\eta-\tilde{\eta})\tilde{\varphi}(\xi,\tilde{\eta})d\tilde{\eta}
-\varsigma(t)\int_{-h_{0}}^{h_{0}}J(\varsigma(t)\eta-\varsigma(t)\tilde{\eta})\tilde{\varphi}(\xi,\tilde{\eta})d\tilde{\eta}\Big)\Big].
\end{array}
\end{align*}
Define
\begin{align*}
\begin{array}{rl}
G(t,\xi,\eta)
=\frac{1+\varepsilon}{\varsigma^{2}(t)}\int_{-h_{0}}^{h_{0}}J(\eta-\tilde{\eta})\tilde{\varphi}(\xi,\tilde{\eta})d\tilde{\eta}
-\varsigma(t)\int_{-h_{0}}^{h_{0}}J(\varsigma(t)\eta-\varsigma(t)\tilde{\eta})\tilde{\varphi}(\xi,\tilde{\eta})d\tilde{\eta}.
\end{array}
\end{align*}
Obviously, $G(t,\xi,\eta)$ is a $T$-periodic function with respect to $\xi$.
Similar as the proof of Theorem 3.3 in \cite{hw19}, we can show that
\begin{align*}
\begin{array}{rl}
G(t,\xi,\eta)
&\geq \frac{\varepsilon}{\varsigma^{2}(t)}\int_{-h_{0}}^{h_{0}}J(\eta-\tilde{\eta})\tilde{\varphi}(\xi,\tilde{\eta})d\tilde{\eta}
-\varsigma(t)\int_{-h_{0}}^{h_{0}}\Big|J(\eta-\tilde{\eta})-J(\varsigma(t)\eta-\varsigma(t)\tilde{\eta})\Big|d\tilde{\eta}\\[3pt]
&\quad-\delta(\delta^{2}+3\delta+3).
\end{array}
\end{align*}
Let
\begin{align*}
\begin{array}{rl}
m=\frac{\varepsilon}{4}\min_{\xi\in[0,T]}\min_{\eta\in [-h_{0}, h_{0}]}\int_{-h_{0}}^{h_{0}}J(\eta-\tilde{\eta})\tilde{\varphi}(\xi,\tilde{\eta})d\tilde{\eta}>0.
\end{array}
\end{align*}
By $(\ref{J})$, there
exists $\delta^{*}\in (0, \frac{1}{2})$ such that for any $0<\delta\leq \delta^{*}$,
\begin{align*}
\begin{array}{rl}
\varsigma(t)\int_{h_{0}}^{h_{0}}\Big|J(\eta-\tilde{\eta})-J(\varsigma(t)\eta-\varsigma(t)\tilde{\eta})\Big|d\tilde{\eta}
\leq \frac{m}{2}.
\end{array}
\end{align*}
It follows that for any $0<\delta\leq \min\{\delta^{*},\frac{m}{10}\}$,
\begin{align*}
G(t,\xi,\eta)\geq 0, \quad\forall (t,\xi,\eta)\in [0,\infty)\times[0,T]\times[-h_{0},h_{0}].
\end{align*}

By the fact $\varsigma(t)\rightarrow 1$ as $\delta\rightarrow 0$, we can choose $0<\sigma, \delta\ll 1$ such that, for
$(t,x)\in[0,\infty)\times (-s(t),s(t))$,
\begin{align*}
\begin{array}{rl}
&\bar{v}_{t}(t,x)-d_{2}[\tau\bar{v}_{xx}+(1-\tau)(\int_{-s(t)}^{s(t)}J(x-y)\bar{v}(t,y)dy-\bar{v}(t,x))]
-\bar{v}(t,x)(c(t)-\bar{v}(t,x))\\[3pt]
&\geq ke^{-\sigma t}\left(-\sigma-\sigma \alpha
+\frac{1}{\varsigma^{2}(t)}\tilde{\lambda}_{1}
+(\frac{1}{\varsigma^{2}(t)}c(\xi)-c(t))\right)\tilde{\varphi}(\xi,\eta)\\[3pt]
&>0.
\end{array}
\end{align*}
Moreover, choosing $k$ large enough such that
\begin{align*}
\bar{v}(0,x)=k\tilde{\varphi}(0,\frac{x}{1+2\delta})\geq v_{0}(x), \quad\forall
x\in [-h_{0},h_{0}].
\end{align*}
Since $s(t)<h_{0}(1+2\delta)<h_{1}$, we know $\bar{u}$ satisfies
\begin{align*}
\begin{array}{rl}
\bar{u}_t\geq d_{1}\int_{-s(t)}^{s(t)}
J(x-y)\bar{u}(t,y)dy-d_{1}\bar{u}+\bar{u}(a(t)-\bar{u}),\quad t>0,
x\in (-s(t), s(t)).
\end{array}
\end{align*}

Note that
\begin{align*}
\begin{array}{rl}
&-\bar{v}_{x}(t,s(t))
=-\frac{k}{\varsigma(t)}e^{-\sigma t}\tilde{\varphi}_{\eta}(\xi(t),h_{0})
\leq \frac{k}{1-\delta}e^{-\sigma t}\|\tilde{\varphi}\|_{C^{1}([0,T]\times [-h_{0},h_{0}])},\\[3pt]
&\int_{-s(t)}^{s(t)}\int_{s(t)}^{\infty}J(x-y)\bar{v}(t,x)dydx
\leq 2kh_{0}(1+2\delta)e^{-\sigma t},\\[3pt]
&\int_{-s(t)}^{s(t)}\int_{s(t)}^{\infty}J(x-y)\bar{u}(t,x)dydx
\leq 2Ch_{0}(1+2\delta)e^{-\frac{\lambda t}{2}}.
\end{array}
\end{align*}
Since $0<\sigma\ll 1$, we may further assume that $\sigma<\frac{\lambda}{2}$.
Suppose that
\begin{align*}
\begin{array}{rl}
0<\mu+\rho_{1}+\rho_{2}\leq \frac{h_{0}\delta\sigma}{2A}
\end{array}
\end{align*}
with
\begin{align*}
\begin{array}{rl}
A:=\max\left\{\frac{k}{1-\delta}\|\tilde{\varphi}\|_{C^{1}([0,T]\times [-h_{0},h_{0}])},
2kh_{0}(1+2\delta), 2Ch_{0}(1+2\delta)\right\},
\end{array}
\end{align*}
we have
\begin{align*}
\begin{array}{rl}
s^{\prime}(t)
&=\frac{1}{2}h_{0}\delta\sigma e^{-\sigma t}\\[3pt]
&\geq \frac{k}{1-\delta}\mu e^{-\sigma t}\|\tilde{\varphi}\|_{C^{1}([0,T]\times [-h_{0},h_{0}])}
+2kh_{0}(1+2\delta)\rho_{2}e^{-\sigma t}+2Ch_{0}(1+2\delta)\rho_{1}e^{-\sigma t}\\[3pt]
&\geq \frac{k}{1-\delta}\mu e^{-\sigma t}\|\tilde{\varphi}\|_{C^{1}([0,T]\times [-h_{0},h_{0}])}
+2kh_{0}(1+2\delta)\rho_{2}e^{-\sigma t}+2Ch_{0}(1+2\delta)\rho_{1}e^{-\frac{\lambda t}{2}}\\[3pt]
&\geq -\mu\bar{v}_{x}(t,s(t))+\rho_{1}\int_{-s(t)}^{s(t)}\int_{s(t)}^{\infty}J(x-y)\bar{u}(t,x)dydx\\[3pt]
&\quad+\rho_{2}\int_{-s(t)}^{s(t)}\int_{s(t)}^{\infty}J(x-y)\bar{v}(t,x)dydx.
\end{array}
\end{align*}
Similarly, we can get
\begin{align*}
\begin{array}{rl}
-s^{\prime}(t)
&\leq -\mu\bar{v}_{x}(t, -s(t))
-\rho_{1}\int_{-s(t)}^{s(t)}\int_{-\infty}^{-s(t)}J(x-y)\bar{u}(t,x)dydx\\
&\quad-\rho_{2}\int_{-s(t)}^{s(t)}\int_{-\infty}^{-s(t)}J(x-y)\bar{v}(t,x)dydx.
\end{array}
\end{align*}
This shows that $(\bar{u}, \bar{v}, -s(t), s(t))$ is an upper solution of $(\ref{1.1})$.
Applying Lemma $\ref{l3.1}$, we have $h(t)\leq s(t)$ and $g(t)\geq -s(t)$, which implies
$h_{\infty}-g_{\infty}\leq 2h_{1}<+\infty$.
\end{proof}

\begin{theorem}\label{t4.12} Assume $a_{T}<d_{1}$. \\
$(i)$ If $h_{0}\geq\frac{1}{2}\min\{h^{*}, l^{*}\}$, then spreading always occurs;\\
$(ii)$ If $h_{0}<\frac{1}{2}\min\{h^{*}, l^{*}\}$,
then there exist $\Lambda^{*}>\Lambda_{*}>0$ such that $h_{\infty}-g_{\infty}<+\infty$ when
$\mu+\rho_{1}+\rho_{2}\leq \Lambda_{*}$ and $h_{\infty}-g_{\infty}=+\infty$ when
$\mu+\rho_{1}+\rho_{2}\geq \Lambda^{*}$.
\end{theorem}

\begin{proof}
$(1)$ If $h_{0}\geq \frac{1}{2}h^{*}$, then Corollary $\ref{r4.9}$ implies that the spreading always occurs.
For the case $h_{0}\geq \frac{1}{2}l^{*}$, if the vanishing happens, then $(g_{\infty},h_{\infty})$
is a finite interval and its length strictly larger than $2h_{0}\geq l^{*}$. Thus,
$\lambda_{1}(-(L_{(g_{\infty},h_{\infty})}+a))<0$, which is a contraction to $(\ref{4.3})$.

$(2)$ From $(\ref{2.2})$, we can deduce that
\begin{align*}
\begin{array}{rl}
&h'(t)>-\mu v_{x}(t, h(t)), \quad
h'(t)>\rho_{1}\int_{g(t)}^{h(t)}\int_{h(t)}^{\infty}J(x-y)u(t,x)dydx, \quad t\geq0,\\[5pt]
&g'(t)<-\mu v_{x}(t, g(t)), \quad
g'(t)<-\rho_{1}\int_{g(t)}^{h(t)}\int_{-\infty}^{g(t)}J(x-y)u(t,x)dydx, \quad t\geq0.
\end{array}
\end{align*}
Since $u,v$ are positive and bounded, we know that $\int_{g(t)}^{h(t)}J(x-y)v(t,y)dy>0$
and there exists a constant $C>0$ such that
$f_{1}\geq -Cu$, $f_{2}\geq -Cv$.
For any given constant $H>\frac{1}{2}\min\{h^{*},l^{*}\}$, by Lemmas 4.1-4.2 in \cite{waw181}, there exist
$\mu^{0}$ and $\rho_{1}^{0}$ such that $g_{\infty}-h_{\infty}\geq 2H$ for any $\mu\geq \mu^{0}$
or $\rho_{1}\geq \rho_{1}^{0}$. Taking $\Lambda^{0}=\mu^{0}+\rho_{1}^{0}$, we know
$g_{\infty}-h_{\infty}=+\infty$ for $\mu+\rho_{1}\geq\Lambda^{0}$.
Applying the continuity method, we can get the desired results.
\end{proof}

Combining Theorems $\ref{t4.8}$, $\ref{t4.10}$-$\ref{t4.12}$ and Corollary $\ref{r4.9}$, we immediately obtain the following criteria for spreading and vanishing.

\begin{corollary}\label{r4.13} (Criteria for spreading and vanishing) Let $(u,v,g,h)$ be the unique solution
of $(\ref{1.1})$, $|\Omega|=h^{*}$ and $|\Omega|=l^{*}$ be the unique root of $\lambda_{1}(-(\tilde{L}_{\Omega}+c))=0$ and
$\lambda_{1}(-(L_{\Omega}+a))=0$, respectively. \\
$(i)$ If one of the following conditions is satisfied:

$(i.1)$ $a_{T}\geq d_{1}$, $(i.2)$ $h_{0}>\frac{1}{2}h^{*}$, $(i.3)$ $a_{T}<d_{1}$ and $h_{0}>\frac{1}{2}l^{*}$,\\
then spreading happens.\\
$(ii)$ If $a_{T}<d_{1}$, $h_{0}<\frac{1}{2}\min\{h^{*},l^{*}\}$,
then there exist $\Lambda^{*}>\Lambda_{*}>0$ such that vanishing happens when
$\mu+\rho_{1}+\rho_{2}\leq \Lambda_{*}$ and spreading happens when
$\mu+\rho_{1}+\rho_{2}\geq \Lambda^{*}$.
\end{corollary}

\section*{Acknowledgments}
Chen¡¯s work was supported by NSFC (No:11801432), China Postdoctoral Science Foundation (No:2019M663610) and the Young Talent fund of University Association for Science and Technology in Shaanxi (No:20200510). Li's work was supported by NSFC (No:11571057). Tang's work was supported by NSFC (Nos:12031010,61772017). Wang's work was supported by NSFC (No:11801429) and the Natural Science Basic Research Plan in Shaanxi Province of China (No:2019JQ-136).

\label{}





\bibliographystyle{model3-num-names}



\end{document}